\documentclass[11pt]{article}
\usepackage[margin=1in]{geometry}
\usepackage{comment}
\usepackage[T1]{fontenc}
\usepackage{lmodern}
\usepackage{todonotes}
\usepackage{ulem}
\usepackage{bm}

\usepackage[colorlinks=true,linkcolor=red,citecolor=blue]{hyperref}%
\usepackage{cite}
\usepackage{url}            
\usepackage{booktabs}       
\usepackage{amsfonts} 
\usepackage{nicefrac}       
\usepackage{microtype} 
\usepackage{wrapfig}
\usepackage{enumerate}
\usepackage{graphicx}
\usepackage{amsmath,amssymb}
\usepackage{bigstrut}
\usepackage{float}
\usepackage{caption}
\usepackage{algorithm}
\usepackage{mathtools}
\usepackage{subfigure}
\usepackage{booktabs}
\usepackage{authblk}
\usepackage{algpseudocode}
\usepackage{dsfont}
\usepackage{bbm}

\newtheorem{theorem}{Theorem}

\newtheorem{remark}{Remark}
\newtheorem{lemma}{Lemma}
\newtheorem{example}{Example}

\newtheorem{proposition}{Proposition}

\newtheorem{asm}{Assumption}

\def\R{\mathbb{R}}
\def\ones{\mathbf{1}}
\newcommand{\EE}[1]{\mathbb{E}\left[ #1 \right]}
\newcommand{\ee}[2]{\mathbb{E}_{#1}\left[ #2 \right]}
\newcommand{\s}{\star}
\newcommand{\bfx}{\mathbf{x}}
\newcommand{\bfy}{\mathbf{y}}
\newcommand{\bfv}{\boldsymbol{v}}

\newcommand{\bfe}{\mathbf{e}}
\newcommand{\bfz}{\mathbf{z}}
\newcommand{\bfr}{\mathbf{r}}

\newcommand{\bfu}{\boldsymbol{u}}

\newcommand{\clg}{\mathcal{G}}

\newcommand{\cle}{\mathcal{E}}

\newcommand{\F}{\mathcal{F}}
\newcommand{\azr}{\alpha_0}
\newcommand{\bzr}{\beta_0}
\newcommand{\kap}{\kappa}
\newcommand{\md}{\middle| }
\DeclareMathOperator{\diag}{diag}

\newcommand{\sign}[1]{\operatorname{sgn}(#1)}

\newcommand{\mx}{\bar{\bfx}}

\newcommand{\lp}{\left(}
\newcommand{\rp}{\right)}
\newcommand{\lb}{\left[}
\newcommand{\rb}{\right]}
\newcommand{\la}{\left\langle}
\newcommand{\ra}{\right\rangle}
\newcommand{\lnr}{\left\|}
\newcommand{\rnr}{\right\|}
\newcommand{\lc}{\left\{}
\newcommand{\rc}{\right\}}
\newcommand{\dl}{\delta}
\newcommand{\et}{\beta(t)}
\newcommand{\at}{\alpha(t)}

\newcommand{\BO}{\mathcal O}

\newcommand{\ER}[1]{E(#1)}
\newcommand{\tu}{\tau}

\def\V{V_\bfr}
\def\rmin{\bfr_{\min}}
\def\Wt{\{W(t)\}}
\def\minW{\eta}

\newcommand{\Del}{\Delta}
\newcommand{\li}{K} 
\newcommand{\sgm}{\sigma}

\newcommand{\OneTwo}{\epsilon_1}
\newcommand{\ThreeOne}{\epsilon_2}
\newcommand{\SumTOne}{\epsilon_3}
\newcommand{\CNeg}{\epsilon_4}

\newcommand{\SumTTwo}{\epsilon_{5}}
\newcommand{\const}{\epsilon_{6}}
\newcommand{\cntB}{\epsilon_{7}(\theta)}
\newcommand{\COne}{\epsilon_{8} }
\newcommand{\CTwo}{\epsilon_{9}} 
\newcommand{\CThree}{\epsilon_{10} }
\newcommand{\CFour}{\epsilon_{11} }
\newcommand{\CBndr}{\epsilon_{12} }
\newcommand{\Cvar}{\epsilon_{13} }

\newcommand{\Nr}[1]{\lnr #1 \rnr_{\bfr}}

\def\R{\mathbb{R}}
\def\ones{\mathbf{1}}
\def\xh{\hat{\bfx}}
\def\nt{\bfe}
\def\sm{\gamma}

\begin{document} 

\title{DIMIX: DIminishing MIXing for Sloppy Agents}

\author{Hadi Reisizadeh, Behrouz Touri, and Soheil Mohajer\footnote{H.\ Reisizadeh (email: hadir@umn.edu) and S.\ Mohajer (email: soheil@umn.edu) are with the University of Minnesota, and B.\ Touri (email: btouri@ucsd.edu) is with the University of California San Diego.}}
\date{}
\maketitle

\begin{abstract}
We study non-convex distributed optimization problems where a set of agents collaboratively solve a separable optimization problem that is distributed over a time-varying network. The existing methods to solve these problems rely on (at most) one time-scale algorithms, where each agent  performs a diminishing or constant step-size gradient descent at the average estimate of the agents in the network. However, if possible at all, exchanging exact information, that is required to evaluate these average estimates, potentially introduces a massive communication overhead. Therefore, a reasonable practical assumption to be made is that agents only receive a rough approximation of the neighboring agents' information. To address this, we introduce and study a \textit{two time-scale} decentralized algorithm with a broad class of \textit{lossy} information sharing methods (that includes noisy, quantized, and/or compressed information sharing) over \textit{time-varying} networks. In our method, one time-scale suppresses the (imperfect) incoming information from the neighboring agents, and one time-scale operates on local cost functions' gradients. We show that with a proper choices for the step-sizes' parameters, the algorithm achieves a convergence rate of $\BO({T}^{-1/3 + \epsilon})$ for non-convex distributed optimization problems over time-varying networks, for any $\epsilon>0$. Our  simulation results support the theoretical results of the paper.
\end{abstract}

\section{Introduction}
Distributed learning serves as a learning framework where a set of computing nodes/agents are interested in collaboratively solving an optimization problem. In this paradigm, in the absence of a central node, the learning task solely depends on on-device computation and local communication among the neighboring agents. With the appearance of modern computation architectures and the decentralized nature of storage, large-scale distributed computation frameworks have received significant attention due to data locality, privacy, data ownership, and scalability to larger datasets and systems. These features of distributed learning have led to applications in several domains including distributed deep networks~\cite{dean2012large, abadi2016tensorflow,jin2016scale}, distributed sensor networks~\cite{li2002detection,rabbat2004distributed,kar2008distributed}, and network resource allocation~\cite{ribeiro2010ergodic,chen2017online}.
In the context of Machine Learning, a prototypical application of distributed learning and optimization is the empirical risk minimization (ERM) problem. In ERM problems, each server in a network of data centers has access to a number of independent data points (realizations) sampled from an underlying distribution, and the objective is to minimize an additive and separable cost function (of the data points) over the underlying network~\cite{nedic2015decentralized,xin2019distributed,zhang2017mixup}.

\noindent\textbf{Related Works.} Decentralized consensus or averaging-based optimization algorithms have been studied extensively over the past few years~\cite{nedic2009distributed,lobel2010distributed,lobel2011distributed,chen2012diffusion, yuan2016convergence,jakovetic2014fast}. It has been shown that when a fixed step-size is utilized, the loss function decreases with the rate of $\BO(1/T)$ until the estimates reach a neighborhood of the (local) minimum of the objective cost function \cite{nedic2009distributed}. However, with a fixed step-size, the local estimates may not converge to an optimal point~\cite{yuan2016convergence}. To remedy this, the diminishing step-size variation of the algorithm is introduced and studied for convex~\cite{nedic2010constrained,jakovetic2014fast,nedic2014distributed} and non-convex problems \cite{sun2016distributed,tatarenko2017non,zeng2018nonconvex}.

A majority of the existing works in this area suffer from requiring large communication overhead, as a (potentially) massive amount of local information is needed to be exchanged among the agents throughout the learning algorithm. In addition, such communication load can lead to a major delay in the convergence of the algorithm,  and become a severe bottleneck when the dimension of the model is large. To mitigate the communication overhead,  several compression techniques are introduced to reduce the size of the exchanged information. In particular, this is used in the context of  distributed averaging/consensus problem~\cite{kashyap2007quantized,nedic2009distributed,cai2011quantized}, where each node has an initial numerical value and aims at evaluating  the average of all initial values using exchange of quantized information, over a fixed or a time-varying network. In the context of distributed optimization, various compression approaches have been introduced to mitigate the communication overhead~\cite{pu2016quantization,reisizadeh2019robust,koloskova2019decentralized,reisizadeh2021distributed,reisizadeh2021adaptive}.

An inexact proximal gradient method with a fixed step-size over a fixed network  for strongly convex loss functions is proposed in~\cite{pu2016quantization}. In this algorithm, each node applies an adaptive deterministic  quantizer on its model and gradient before sending them to its neighbors. It is shown that under certain conditions, 
the algorithm provides a linear convergence rate. In another  related work~\cite{reisizadeh2019robust}, authors proposed a decentralized gradient descent method, named `QuanTimed-DSGD', with fixed step-sizes. In this algorithm, agents exchange a quantized version of their models/estimates in order to reduce the communication load. 
It is shown that for strongly convex loss functions and \textit{a given termination time} $T$, constant step-size parameters for the algorithm can be chosen  so that the average distance between the model of the agents and the global optimum is bounded by $c(T^{-1/2+\epsilon})$ for some  $c>0$ and any  $\epsilon>0$. In this algorithm not only the step sizes depend on $T$, but also the results only hold for the exact termination iteration $T$, which further needs to satisfy $T\geq T_{\min}$, where $T_{\min}$ itself depends on $\epsilon$ as well as  non-local parameters of the \textit{fixed} underlying network. 
For smooth non-convex cost functions, there it is shown that the temporal average (up to time $T$) of the expected norm of the gradients at average models does not exceed $c(T^{-1/3})$. However, similar to the strongly convex case, the choice of parameters (i.e., the step-sizes) for QuanTimed-DSGD depend on the termination time $T$, and need to be re-evaluated when the target termination time varies. To compensate the quantization error, a decentralized (diminishing) gradient descent algorithm is proposed in  in~\cite{koloskova2019decentralized,koloskova2019decentralizedICLR} using error-feedback. The proposed algorithm achieve the convergence rate of $\BO(T^{-1})$ and $\BO(T^{-1/2})$ for strongly and non convex objective functions, respectively. However, the nature of the algorithm restricts its use to time-invariant networks and in addition, the feedback mechanism cannot compensate communication noise between the nodes. In~\cite{srivastava2011distributed}, a two-time-scale gradient descent algorithm was presented for distributed constrained and convex optimization problems over an i.i.d.\  communication graph with noisy communication links, and sub-gradient errors. It is shown that if the random communication satisfies certain conditions,  proper choices of the time-scale parameters result in the almost sure convergence of the local states to an optimal point. Another interesting approach to address exact convergence for distributed optimization with  fixed gradient step-sizes under a \textit{noiseless} communication model is to use gradient tracking methods~\cite{pu2018distributed,zhang2019decentralized}.

\noindent\textbf{Our Contributions.}  
We introduce and study a \textit{two time-scale} decentralized gradient descent algorithm for a \textit{broad} class of imperfect sharing of information over \textit{time-varying communication networks} for distributed optimization problems with smooth \textit{non-convex} local cost functions. In this work, one time-scale addresses the convergence of the local estimates to an stationary point while the (newly introduced) second time-scale is introduced to suppress the noisy effect of the imperfect incoming information from the neighbors. 

Our main result shows that with a proper choice of the parameters for the two diminishing step-size sequences, the temporal average of the expected  norm of the gradients decreases with the rate of $\BO(T^{-1/3+\epsilon})$. To prove this result, we have established new techniques to analyze  the interplay {between} the two time-scales, in particular, in the presence of underlying time-varying networks. To validate our theoretical results, we numerically evaluate the algorithm for optimizing a regulized logistic regression on benchmark dataset MNIST. To assess the efficiency of our algorithm over time-varying networks, we also implement a linear regression model (convex) over synthetic data and a CNN (non-convex) over CIFAR-10. Our simulation results also confirm that introducing a \textit{time-shift} in our step-sizes improves the \textit{finite-time} performance of our algorithm leading to a faster decay of the loss function.

\noindent\textbf{Paper Organization.} 
After introducing some notations, we present the problem formulation, the algorithm, the main result, and some of its practical implications in Section~\ref{sec:prob-def-main-res}. Our theoretical results are corroborated by simulation results in Section~\ref{sec:expr-res}. To prove the main result, we first provide some auxiliary lemmas in Section~\ref{sec:aux-lemmas} whose proofs are presented in Appendix~\ref{sec:proof-aux}. In~Sections~\ref{sec:proof_nn_cvx}, we present the proof of the main result. Finally, we conclude the paper in Section~\ref{sec:conclusions}.

\noindent\textbf{Notation and Basic Terminology.} 
Throughout this paper, we use $[n]$ to denote the set of integers $\{1,2,\dots,n\}$.  In this paper, we are dealing with $n$ agents that are minimizing a function in $\R^d$. For notational convenience, throughout this paper, we assume that the underlying functions are acting on \textbf{row} vectors, and hence, we view vectors in  $\R^{1\times d}=\R^d$ as row vectors. The rest of the vectors, i.e., the vectors in $\R^{n\times 1} = \R^n$, are assumed to be column vectors.   The $\ell_2$-norm of  a vector $\bfx\in \mathbb{R}^d$ is defined as $\|\bfx\|^2=\sum_{j=1}^d |x_j|^2$. The Frobenius norm of a matrix $A\in \mathbb{R}^{n\times d }$ is defined as  $\lnr A\rnr_F^2 = \sum_{i=1}^n \|A_i\|^2 =\sum_{i=1}^n \sum_{j=1}^d |A_{ij}|^2$. A vector $\bfr\in \mathbb{R}^n$ is called stochastic if $r_i\geq 0$ and $\sum_{i=1}^n r_i=1$. Similarly, a non-negative matrix $A\in \mathbb{R}^{n\times d}$ is called (row) stochastic if $\sum_{j=1}^d A_{ij}=1$ for every $i\in [n]$. 
For a matrix $A\in \mathbb{R}^{n\times d}$, we denote its $i$-th row and $j$-th column by $A_i$ and $A^j$, respectively. Note that we deal with two types of vector throughout the paper.  For an $n\times d$  matrix $A$ and a strictly positive stochastic vector $\bfr\in \mathbb{R}^n$, we define the $\bfr$-norm of $A$  by $\lnr A \rnr_{\bfr}^2 = \sum_{i=1}^n r_i \lnr A_i\rnr^2$. It can be verified that $\Nr{\cdot}$ is a norm on the space of $n\times d$ matrices.  Finally, we write $A \geq B$ if $A-B$ (is well-defined and) has non-negative entries. 
\section{Problem Setup and Main Result}\label{sec:prob-def-main-res}
In this section, first we formulate non-convex distributed optimization problems over time-varying networks and introduce some standard assumptions on the underlying problem. After proposing  our algorithm, we  state our  main result. Finally, we discuss the implications of our result on various important practical settings with imperfect information sharing.  

\subsection{Problem Setup}
This paper is motivated by stochastic learning problems in which the goal is to solve
    \begin{align}\label{eq:st_l_f}
    \min_{\bfx} L(\bfx) := \min_{\bfx} \ee{\xi \sim \mathcal{P}}{\ell(\bfx,\xi)},
    \end{align}
    where $\ell:\mathbb{R}^d\times \mathbb{R}^p\rightarrow \mathbb{R}$ is a loss function, $\bfx\in \mathbb{R}^{1\times d} = \R^d$ is the decision/optimization row vector,  and $\xi$ is a random vector taking values in $\mathbb{R}^p$ that is drawn from an unknown underlying distribution $\mathcal{P}$. One of the key practical considerations that renders~\eqref{eq:st_l_f} as a challenging task is that the underlying distribution $\mathcal{P}$ is often unknown. Instead, we have access to~$N$ independent realizations of $\xi$ and focus on  solving the corresponding ERM problem which is given by
    \begin{align}\label{eq:ERM}
    \min _{\bfx} f(\bfx) := \min_{\bfx} \frac{1}{N}\sum_{j=1}^{N}\ell(\bfx,\xi_j),
    \end{align}
    where $f(\bfx)$ is the empirical risk with respect to the data points $\mathcal{D}=\{\xi_1,\ldots,\xi_N\}$. We assume that $\ell(\cdot,\cdot)$ is a non-convex loss function, which potentially results in a non-convex function $f(\cdot)$. 
   
    In distributed  optimization, we have a network consisting of $n$ computing nodes (agents, or workers), where each node~$i$  observes a non-overlapping subset of ${m_i = r_i N}$  data points, denoted by $\mathcal{D}_i = \{\xi^i_1,\ldots,\xi^i_{m_i}\}$, where  ${\mathcal{D} = \mathcal{D}_1\cup\cdots\cup \mathcal{D}_n}$. Here, $r_i$ represents the fraction of the data that is processed at node $i\in[n]$. Note that the vector $\bfr=(r_1,\dots, r_n)$ is a strictly positive stochastic vector, i.e., $r_i>0$ and $\sum_{i=1}^n r_i=1$. Thus, the ERM problem in~\eqref{eq:ERM} can be written as the minimization of the weighted average of local empirical risk functions $f_i$ for all nodes $i\in [n]$ in the network, i.e.,
    \begin{align}\label{eq:ERM2}
    \min_{\bfx}  f(\bfx) =\min_{\bfx} \sum_{i=1}^{n}r_i f_i(\bfx)=\min_{\bfx}\frac{1}{N}\sum_{i=1}^{n}\sum_{\xi\in\mathcal{D}_i}\ell(\bfx,\xi),
    \end{align}
     where $f_i(\mathbf{x}) :=\frac{1}{m_i}\sum_{\xi\in\mathcal{D}_i}\ell(\bfx,\xi)= \frac{1}{m_i}\sum_{j=1}^{m_i}\ell(\mathbf{x},\xi_{j}^{i})$.
	We can rewrite the ERM problem in~\eqref{eq:ERM2} as a distributed consensus optimization problem, given by
	\begin{align}\label{eq:ERM3}
	 \min_{\bfx_1,\ldots,\bfx_n} \sum_{i=1}^{n}r_i f_i(\bfx_i) \quad \textrm{subject to}\quad \bfx_1=\bfx_2=\cdots=\bfx_n.
	\end{align}
	Consider an $n\geq 2$ agents that are connected through a \textit{time-varying} network. We represent this network at time $t\geq 1$ by the directed graph  $\clg(t)=([n], \cle(t))$, where the vertex set $[n]$ represents the set of agents and the edge set ${\cle(t) \subseteq \{(i,j):i,j\in [n]\}}$ represents the set of links at time $t$. At each discrete time $t\geq 1$, agent~$i$ can only send information to its (out-) neighbors in $\cle(t)$, i.e., all $j\in [n]$ with $(i,j) \in \cle(t)$.

    To discuss our algorithm (\texttt{DIMIX}) for solving \eqref{eq:ERM3} distributively, let us first discuss its general structure and the required information at each node for its execution. In this algorithm, at each iteration $t\geq 1$, agent $i\in [n]$ updates its estimate $\bfx_i(t)\in \R^d$ of an optimizer of \eqref{eq:ERM2}. To this end, it utilizes the  gradient information of its own local cost function $f_i(\bfx)$ as well as a \textit{noisy/lossy} average of its current neighbors estimates, denoted by 
    ${\hat{\bfx}_i(t) := \sum_{j=1}^n W_{ij}(t)\bfx_j(t)+\nt_i(t)}$. Here,  $W(t)$ is a \textit{row-stochastic} matrix that is consistent with the underlying connectivity network~$\clg(t)$ (i.e., $W_{ij}(t)>0$ only if $(j,i)\in \cle(t)$) and $\nt_i(t)\in \R^d$ is a random {noise} vector. Later, {in Section~\ref{sub:examples}} we discuss several noisy and lossy information sharing architectures (quantization and noisy
    communication) that fit in this broad information structure. Now we are ready to discuss the \texttt{DIMIX} algorithm. In this algorithm, using the information available to agent $i$ at time $t$, agent $i$ updates its current estimate by computing  a \textit{diminishing} weighted average of its own state and the noisy average of its neighbors' estimates, and moves along its local gradient. More formally, the update rule at node $i\in [n]$ is given by
   \begin{align}\label{eq:upd_nd}
    \bfx_i(t+1) = (1-\beta(t))\bfx_i(t)+ \beta(t)\xh_i(t)- \alpha(t) \beta(t)\nabla f_i(\bfx_i(t)),
   \end{align}
   where $\alpha(t) = \frac{\azr}{(t+\tu)^{\nu}}$ and $\beta(t) = \frac{\bzr}{(t+\tu)^\mu}$ for some  $\mu, \nu \in (0,1)$ are the diminishing step-sizes of the algorithm, and $\tu\geq 0$ is an arbitrary shift, that is introduced to accelerate the finite-time performance of the algorithm. The description of \texttt{DIMIX} is summarized in Algorithm~\ref{alg:1}.
   For notational simplicity, let
       \begin{align}
        X(t) :=\begin{bmatrix} 
        \bfx_1(t) \\ \vdots \\ \bfx_n(t)
        \end{bmatrix},\quad 
        E(t) := \begin{bmatrix}
        \nt_1(t) \\ \vdots\\ \nt_n(t)\end{bmatrix},\quad 
        \nabla f(X(t)) := \begin{bmatrix} \nabla f_1(\bfx_1(t))\\ \vdots \\ \nabla f_n(\bfx_n(t))\end{bmatrix}.
        \label{eq:def:grad-f-X}
    \end{align}
    Since $\hat{\bfx}_i(x) = \sum_{j=1}^n W_{ij}(t) \bfx_j(t) + \bfe_i(t)$, we can rewrite the update rule in~\eqref{eq:upd_nd} in the matrix format
    \begin{align}
        X(t+1) = ((1-\et)I+\et W(t))X(t)+\et\ER{t} - \at\et\nabla f(X(t)). 
        \label{eq:dyn-col-mat}
    \end{align}
   Let us discuss  some important aspects of the above update rule. Note that both $\alpha(t)$ and $\beta(t)$ are diminishing step-sizes. If $\beta(t)=\beta_0<1$ and $\alpha(t)=\alpha_0<1$ are both constants, then the dynamics in~\eqref{eq:dyn-col-mat} reduces to the algorithm proposed in~\cite{reisizadeh2019robust} for both convex and non-convex cost functions. Alternatively, if $\beta(t)=\beta_0<1$ is a constant sequence and $E(t)=\mathbf{0}$ for all $t\geq 1$, \eqref{eq:dyn-col-mat} would be reduced to the averaging-based distributed optimization with diminishing steps-sizes (for gradients), which is introduced and studied in \cite{nedic2010constrained} for local convex cost functions $f_i(\bfx)$.  The newly introduced time-scale/step-size $\beta(t)$ suppresses the incoming noise $\bfe_i(t)$ from the neighboring agents. However, $\beta(t)$ also suppresses  the incoming signal level  $\sum_{j=1}^nW_{ij}(t)\bfx_j(t)$ at each node~$i$. This casts a major technical challenge for establishing convergence-to-consensus guarantees for the algorithm over time-varying networks. 
   
   On the other hand, the diminishing step-size for the gradient update is $\hat{\alpha}(t)=\alpha(t)\beta(t)$. We chose to represent our algorithm in this way to ensure that the local mixing (consensus) scheme is operated on a faster time-scale than the gradient descent. 
   
    \begin{algorithm}[t!]
    \caption{\texttt{DIMIX} at agent $i$}\label{alg:1}
    \begin{algorithmic}[1]
    \Require Stochastic matrix sequence $\{W(t)\}$, Iteration $T$
    \State Set $\bfx_{i}(1)=0$. 
    \For{$t=1,\ldots,T-1$}
    \State Compute the local gradient $\nabla f_i(\bfx_{i}(t))$.
    \State Obtain noisy average neighbors' estimate $\xh_i(t)$.
    \State Set: 
      $ \bfx_{i}(t+1) = (1-\beta(t))\bfx_i(t)+ \beta(t)\xh_i(t)- \alpha(t) \beta(t)\nabla f_i(\bfx_i(t))$.
\EndFor
\end{algorithmic}
\end{algorithm}
\subsection{Assumptions}

In order to provide performance guarantees for \texttt{DIMIX}, we need to assume certain regularity conditions on (i)~the statistics of the (neighbors' averaging) noise process $\{E(t)\}$,  (ii)~the mixing properties of the weight sequence $\{W(t)\}$, and (iii)~the loss function $\ell(\cdot,\cdot)$.

First, we discuss our main assumption on the noise sequence $\{E(t)\}$. 
\begin{asm}[\textbf{Neighbors State Estimate Assumption}]\label{assum:neighbor}
Suppose that $\{X(t)\}$ is adapted to a filtration $\{\F_t\}$ on the underlying probability space (see e.g.\ Section 5.2 in \cite{durrett2019probability}). {We assume that there exists some $\sm>0$ such that for all $i\in [n]$ and all $t\geq 1$, the noise sequence $\{\nt_i(t)\}$ satisfies} 
        \begin{align}\label{eqn:condexp}
            \EE{\nt_i(t)\mid \F_t}&=0\mbox{, and}\cr 
            \EE{\|\nt_i(t)\|^2 \mid \F_t}&\leq \sm.
        \end{align}
 \end{asm}
 Note that the natural filtration of the random process $\{X(t)\}$ is one choice for $\{\F_t\}$. In this case,~\eqref{eqn:condexp} reduces to $\EE{\nt_i(t)\!\mid\!\! X(1),\ldots,X(t)}\!=\!0$ and  $\EE{\|\nt_i(t)\|^2 \!\mid\! X(1),\ldots,X(t)}\!\leq \sm.$

 Next, we discuss the main assumption on the network connectivity, or more precisely, the mixing of information over the time-varying network. 
\begin{asm}[\textbf{Connectivity Assumption}]\label{asm:W}
 We assume that the weight matrix sequence $\Wt$ in \eqref{eq:dyn-col-mat} satisfies the following properties
 \begin{enumerate}[(a)]
     \item \textit{Stochastic with Common Stationary Distribution}: $W(t)$ is non-negative and  ${W(t)\ones=\ones}$ and ${\bfr^T W(t)=\bfr^T}$ for all $t\geq 1$, where $\ones\in \R^n$ is the all-one vector, and $\bfr>0$ is the weight vector. 
     
     \item \textit{Bounded Nonzero Elements}: There exists some $\minW>0$ such that if for some $i,j\in [n]$ and $t\geq 1$ we have  $W_{ij}(t)>0$, then $W_{ij}(t)\geq \minW$.
     
     \item \textit{$B$-Connected}: For a fixed integer $B\geq 1$, the graph $\lp [n],\bigcup_{k=t+1}^{t+B}\cle(k)\rp$ is strongly connected for all $t\geq 1$, where $\cle(k)=\{(j,i)\mid W_{ij}(k)>0\}$. 
 \end{enumerate}
 \end{asm}
 
 Finally for the loss function $\ell(\cdot,\cdot)$, we make the following assumption. 
 
\begin{asm}[\textbf{Stochastic Loss Function Assumption}]\label{asm:f}
We assume that:
\begin{enumerate}[(a)]
      \item	The function $\ell(\cdot, \cdot)$ is $\li$-smooth with respect to its first argument, i.e., for any $\bfx,\bfy\in\mathbb{R}^d$  and $\xi\in \mathcal{D}$ we have that $\|\nabla\ell(\bfx,\xi)-\nabla\ell(\bfy,\xi)\|\leq \li\|\bfx-\bfy\|$.
      \item Stochastic gradient $\nabla\ell(\bfx,\xi)$ is unbiased and variance bounded, i.e.,
      \[\ee{\xi}{\nabla\ell(\bfx,\xi)} = \nabla L(\bfx), \quad \ee{\xi}{\|\nabla\ell(\bfx,\xi) - \nabla L(\bfx)\|^2}\leq \sgm^2.\]
\end{enumerate}
\end{asm}
Note that Assumption~\ref{asm:f}-b implies a homogeneous sampling, i.e., each agent draws i.i.d.\ sample points from a data batch.  In a related work~\cite{lian2017can}, a stronger assumption has been considered which allows for heterogeneous data samples.
\begin{remark}
Recall that local function $f_i$ is defined as $f_i(\bfx) = \frac{1}{m_i} \sum_{j=1}^{m_i} \ell(\bfx, \xi_{j}^{i})$. Then the conditions of Assumption~\ref{asm:f} on the function $\ell(\cdot,\cdot)$ can be directly translated to conditions on $f_i(\bfx)$ given by
\begin{align*}
    \| f_i(\bfx) - f_i(\bfy) \| \leq \li \|\bfx - \bfy\|,\quad
    \ee{\xi}{\nabla f_i(\bfx)} = \nabla L(\bfx),\quad
    \ee{\xi}{\|\nabla f_i(\bfx) - \nabla L(\bfx)\|^2}\leq \sgm^2/m_i.
\end{align*}.
\end{remark}
\subsection{Main Result}
Here, we characterize the convergence rates of our algorithm for the $K$-smooth non-convex loss functions. More precisely, we establish a rate for the temporal average of the expected norm of the gradients for various choice of the time-scale parameters $\nu,\mu$. 
\begin{theorem}\label{thm:non_cvx}
Suppose that Assumptions~\ref{assum:neighbor}--\ref{asm:f} hold and let $\alpha(t) = \frac{\azr}{(t+\tu)^{\nu}}$, $\beta(t) = \frac{\bzr}{(t+\tu)^\mu}$ where $\alpha_0, \beta_0\in (0,1)$, $\tau\geq 0$, and $\nu,\mu \in (0,1)$ are arbitrary constants with $\mu\neq 1/2$ and $3\nu+\mu \neq 1$. Then 
the weighted average estimates $\mx(t):= \sum_{i=1}^{n}r_i\bfx_i(t)$ generated by Algorithm~\ref{alg:1} satisfy 
\begin{align}\label{eq:thm:rate}
   M_\theta(\nu,\mu):= \lb \frac{1}{T}\sum_{t=1}^{T} \left( \EE{\lnr\nabla f(\mx(t))\rnr^2} \right)^\theta\rb^{1/\theta}
   &= \BO\lp T^{-\min\{1-\nu-\mu, \mu-\nu, 2\nu\} }\rp,
\end{align}
where $\theta\in (0,1)$ is an arbitrary constant. 

Furthermore, for $(\nu^\star, \mu^\star)=\lp \frac{1}{6}, \frac{1}{2}\rp$  we get the optimal rate of
\begin{align}\label{eq:thm:rate:opt}
   M_\theta (\nu^\star, \mu^\star) = \BO\lp T^{-1/3}\ln T\rp.
\end{align}
\end{theorem}
\begin{remark}\label{rem:randomness}
Note that the expectation operator $\EE{\cdot}$ is over the randomness of the dataset $\mathcal{D}$ and the compression/communication noise. Moreover, note that the theorem above shows that the gradient of $f(\cdot)$ (which depends of the choice of $\mathcal{D}$) at the average state of $\overline{\mathbf{x}}(t)$ (which also depends on $\mathcal{D}$) vanishes at a certain rate. It is worth mentioning that  this is not the performance of the average trajectory for the average function. 
\end{remark}
\begin{remark}
From~\eqref{eq:thm:rate}, one has to maximize $\min\{1-\nu-\mu, \mu-\nu, 2\nu\}$ over $\nu, \mu\in(0,1)$  to achieve the fastest convergence for $M_\theta$. This leads to $(\nu^\star, \mu^\star)=(1/6,1/2)$, which none of the conditions $\mu\not= 1/2$ and $3\nu+\mu \not= 1$ hold for. However, one can choose $(\nu, \mu)=(1/6+\epsilon/2,1/2+\epsilon/2)$  and obtain $M_{\theta} = \BO(T^{-1/3+\epsilon})$ for any $\epsilon>0$. Nevertheless, note that \eqref{eq:thm:rate:opt} provides a faster convergence rate of $\BO(T^{-1/3} \ln T)$ for $(\nu^\star,\mu^\star)=(1/6,1/3)$. 
\end{remark}
\begin{proposition}\label{cor:M1}
Under the conditions of Theorem~\ref{thm:non_cvx}, for the optimum choice of $(\nu^\star,\mu^\star)=(1/6,1/3)$, we have
\begin{align}\label{eq:thm:M1}
   M_1 (\nu^\star,\mu^\star) :=  \frac{1}{T}\sum_{t=1}^{T}  \EE{\lnr\nabla f(\mx(t))\rnr^2}  
   &\leq \BO\lp T^{-1/3+\epsilon}\rp,
\end{align}
for any $\epsilon>0$.  
Furthermore, in this case, for each agent $i\in [n]$ the convergence rate to consensus is given by 
\begin{align}\label{eq:thm:var}
    \frac{1}{T}\sum_{t=1}^{T}\EE{\|\bfx_i(t)-\bar{\bfx}(t)\|^2} \leq \BO\lp T^{-1/3+\epsilon}\rp.
\end{align}
As a result, combining~\eqref{eq:thm:M1},~\eqref{eq:thm:var}, and Assumption~\ref{asm:f}, for all $i\in [n]$, we have
\begin{align*}
    \frac{1}{T}\sum_{t=1}^{T}  \EE{\lnr\nabla f(\bfx_i(t))\rnr^2}  
   &\leq \BO\lp T^{-1/3+\epsilon}\rp.
\end{align*}
\end{proposition}

\begin{remark}\label{rem:vanishing-graph}
We should comment that almost universally, all the existing results and algorithms on distributed optimization algorithms for time-varying graphs assume a uniform positive lower bound on non-zero elements of the (effective) weight matrices \cite{nedic2008distributed,nedic2009distributed,nedic2013distributed,touri2013product,touri2014endogenous,yuan2016convergence,tatarenko2017non,aghajan2020distributed}. Absence of such an assumption significantly increases the complexity of the convergence analysis of the algorithm. In our work, even though the stochastic matrix sequence $\{W(t)\}$ is assumed to be $B$-connected, the \textbf{effective averaging} weight sequence, that is given by $\{(1-\beta(t))I+\beta(t)W(t)\}$,  has vanishing (non-self) weights. One of the major theoretical contributions of this work is to introduce tools and techniques to study distributed optimization with diminishing weight sequences.
\end{remark}

\begin{remark}
In a related work~\cite{reisizadeh2019robust} on distributed optimization with compressed  information sharing among the nodes, authors considered a fixed step-size (zero time-scale) version of Algorithm~\ref{alg:1} with a fixed averaging matrix $W$. It is shown that for a given \textbf{termination time} $T$, the algorithm's step-sizes can be chosen (depending on $T$)  such that the temporal average (up to iteration $T$) of the expected norm of the gradient (i.e., $M_1$ defined in~\eqref{eq:thm:M1}) does not exceed $c({T}^{-1/3 })$ (where $c>0$ is a constant). However, the step-sizes should be re-evaluated and the algorithm needs to be re-executed  if one targets another termination time $T'$. In this work, we use vanishing step-sizes $\alpha(t)$ and $\beta(t)$ (which do not depend on the termination time) and show that the same temporal average vanishes at the rate of $\BO\lp{T}^{-1/3+\epsilon}\rp$ for \textbf{every} iteration $T$ and any arbitrarily small  $\epsilon>0$.
\end{remark}
\subsection{Examples for Stochastic Noisy State Estimation}\label{sub:examples}

\begin{figure}
    \centering
    \includegraphics[width=\textwidth]{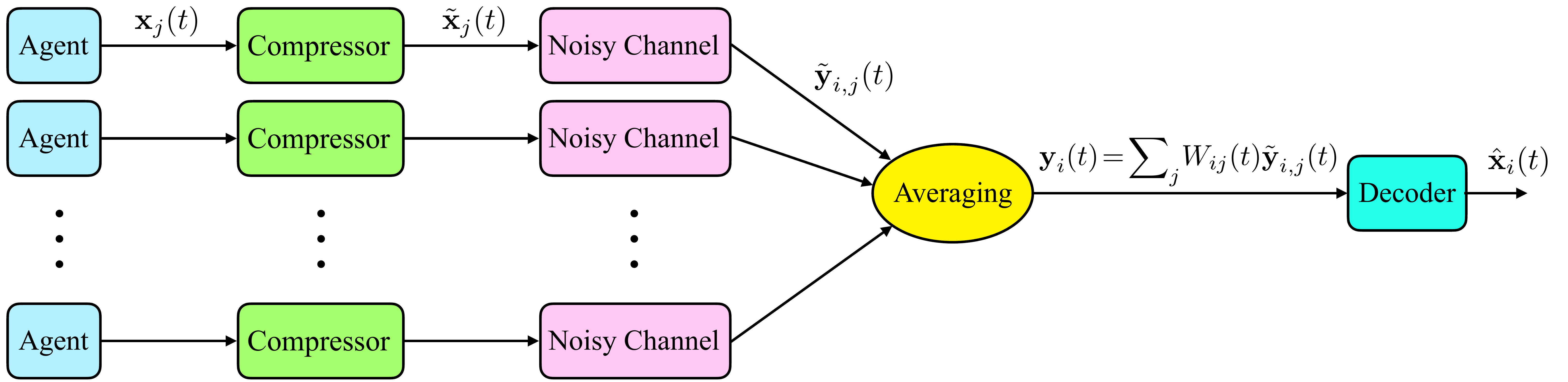}
    \caption{A general architecture to lossy information model.}
    \label{fig:architecture}
\end{figure}
The noisy information in~\eqref{eq:upd_nd} is very general, and captures several models of imperfect data used in practice and/or theoretical studies. A rather general architecture that leads to such noisy/lossy estimates is demonstrated in Figure~\ref{fig:architecture}: Once  the estimate $\bfx_j(t)$ of node $j$ at iteration $t$ is evaluated,  node $j$ may apply an operation (such as compression/sparsification or quantization) on its own model to generate $\tilde{\bfx}_j(t)$. This data will be sent over a potentially noisy communication channel, and a neighbor $i$ will receive a corrupted version of $\tilde{\bfx}_j(t)$, say $\tilde{\bfy}_{i,j}(t)$ from every neighbor node $j$.  Upon receiving such channel outputs from all its neighbors, node $i$ computes their weighted average $\bfy_i(t) = \sum_{j=1}^n W_{ij}(t) \tilde{\bfy}_{i,j}(t)$. This  can be decoded to the approximate average model $\hat{\bfx}_i(t)$. In the following we describe three popular frameworks, in which each node~$i$ can only use an imperfect neighbors' average $\hat{\bfx}_i(t)$ to update its estimate. 
It is worth emphasizing that these are just some examples that lie under the general model in~\eqref{eq:upd_nd}.

\begin{example}\textbf{(Unbiased Random Sparsification).}
Let $Z_k$ be the set of all subsets of $[d]$ of size $k$ for a  parameter $k\in [d]$. 
We define the operator  ${\text{Rand}_k:\R^d\times Z_k\rightarrow\R^d}$ as 
\begin{align}
    [\text{Rand}_k(\bfx,\zeta)]_\ell :=\left\{
    \begin{array}{ll}
       (\frac{d}{k}) x_{\ell} & \ell\in \zeta\\
       0 & \textrm{otherwise}.
    \end{array}\right.
\end{align}
 For a given time, with $\zeta_j$ being drawn uniformly at random from $Z_k$ independent of the input vector $\bfx$, iteration,  and other agents, $\text{Rand}_k(\cdot,\cdot)$  is an unbiased version of the random sparisification introduced in~\cite{stich2018sparsified}. 
 It can be verified that, when  ${\|\bfx\|^2\leq D}$, we have
 \[\mathbb{E}_\zeta\lb \text{Rand}_k(\bfx,\zeta)\rb=\bfx,\quad \mathbb{E}_\zeta\lb\|\text{Rand}_k(\bfx,\zeta)-\bfx\|^2\rb\leq \lp \frac{d}{k}-1 \rp \|\bfx\|^2 \leq \lp \frac{d}{k} -1 \rp D.\] 
 Then, we can get
    	\begin{align*}
    	    \hat{\bfx}_i(t) = \sum_{j=1}^n W_{ij}(t) \text{Rand}_k(\bfx_j(t),\zeta_j)  = \sum_{j=1}^n W_{ij}(t) \bfx_j(t) + \bfe_i(t),
    	\end{align*}
    	where  
    	$\nt_i(t) = \sum_{j=1}^n W_{ij}(t) \lp \text{Rand}_k(\bfx_j(t),\zeta_j) - \bfx_j(t)\rp$. From the properties of $\text{Rand}_k(\cdot,\cdot)$, we can verify that $\EE{\nt_i(t)|\F_t} =0$ and 
    	\[
    	\EE{\|\nt_i(t)\|^2|\F_t} = \lp \frac{d}{k}-1 \rp D \sum_{j=1}^n W_{ij}(t)^2 \leq \lp \frac{d}{k}-1 \rp D,
    	\]
    	which fulfills the conditions of Assumption~\ref{assum:neighbor}. With respect to the architecture of Figure~\ref{fig:architecture}, we have $\tilde{\bfx}_j(t) = \text{Rand}(\bfx_j(t),\zeta_j)$. Moreover, the noisy channel is perfect, and the decoder is just an identity function, i.e., $\bfy_{i,j}(t)=\tilde{\bfx}_j(t)$ and $\hat{\bfx}_i(t) = \bfy_i(t)$. 
\end{example}
\begin{example}\label{exm:SQb}\textbf{(Stochastic Quantizer with bounded trajectory).} 
    	The stochastic quantizer with a number of quantization levels $s$ maps a vector $\bfx\in \mathbb{R}^d$ to a random vector  $Q_s^{S}(\bfx) \in \R^d$,  where its $\ell$th entry is given by
    	\begin{align}\label{eq:unb-q}
    	\left[Q^{S}_s(\bfx)\right]_\ell := \|\bfx\|\cdot \sign{x_\ell}\cdot \zeta\lp |x_\ell|/\|\bfx\|,s \rp ,\quad  \ell\in[d],
    	\end{align}
    	and $\zeta(x, s)$ is a random variable taking values 
        \begin{align}\label{eq:zeta}
    	    \zeta(x,s) = \left\{
    	    \begin{array}{ll}
    	    \lceil s x \rceil / s & \textrm{w.p. $s x - \lfloor s x\rfloor$}\\
    	    \lfloor s x \rfloor / s & \textrm{w.p. $\lceil s x\rceil - s x$.}
    	   \end{array}\right.
    	\end{align}
    	Note that, random variables $\{\zeta(\cdot,\cdot)\}$ are independent, across the coordinates, agents, and time steps. Thus, in this case, the relationship between $\tilde{\bfx}_j(t)$ and $\bfx_j(t)$ in Figure~\ref{fig:architecture} would be $\tilde{\bfx}_j(t)= Q^{S}_s(\bfx_j(t))$. Furthermore, the noisy channel is perfect, and the decoder component is just an identity function, i.e., $\bfy_{i,j}(t)=\tilde{\bfx}_j(t)$ and $\hat{\bfx}_i(t) = \bfy_i(t)$. 
    	It is shown in~\cite{alistarh2017qsgd} that applying this quantizer on $\bfx \in \mathbb{R}^d$ with a  bounded norm $\|\bfx\|^2\leq D$ satisfies $\EE{Q^{S}_s(\bfx)}=\bfx$ and $\EE{\|Q^{S}_s(\bfx) - \bfx\|^2}\leq \min\left(\frac{\sqrt{d}}{s}, \frac{d}{s^{2}}\right)D$.
    	Therefore, the neighbors estimate for node $i$ will be
    	\begin{align*}
    	    \hat{\bfx}_i(t) = \sum_{j=1}^n W_{ij}(t) Q^{S}_s(\bfx_j(t))  = \sum_{j=1}^n W_{ij}(t) \bfx_j(t) + \bfe_i(t),
    	\end{align*}
    	where $\nt_i(t) \!= \!\sum_{j=1}^n W_{ij}(t) \lp Q^{S}_b(\bfx_j(t)) \!-\! \bfx_j(t)\rp$ satisfies
    	$\EE{\nt_i(t)|\F_t} =0$ and 
    	\[
    	\EE{\|\nt_i(t)\|^2|\F_t} = \min\left(\sqrt{d}/s, d/s^{2}\right)D \sum_{j=1}^n W^2_{ij}(t) \leq \min\left(\sqrt{d}/s, d/s^{2}\right)D,
    	\]
    	provided that $\|\bfx_j(t)\|^2\leq D$ for every $j\in[n]$ and every $t\geq 1$. Therefore, the conditions of Assumption~\ref{assum:neighbor} are satisfied.  
	\end{example}
	Note that in both examples above $\nt_i(t)$ and $\nt_j(t)$ might be correlated, especially when nodes~$i$ and~$j$ have common neighbor(s). However, this does not violate the conditions of Assumption~\ref{asm:f}.

\begin{example}
	\textbf{(Noisy Communication).}
	The noisy neighbor estimate model may arise due to imperfect communication between the agents.  Consider a wireless network, in which the computing nodes communicate with their neighbors over a Gaussian channel, i.e., when node $j$ sends its state $\bfx_j(t)$ (without  compression, i.e., $\tilde{\bfx}_j(t) = \bfx_j(t)$) to its neighbor $i$,  the signal received at node $i$ is $\tilde{\bfy}_{i,j}(t)=\bfx_j(t) + \bfz_{i,j}(t)$ where $\bfz_{i,j}(t)$ is a zero-mean Gaussian noise with variance $\zeta^2$, independent across $(i,j)$, and $t$. 
	Applying an identity map decoder at node~$i$ (i.e., $\hat{\bfx}_i(t) = \bfy_i(t)$) we have 
	\begin{align*}
	    \hat{\bfx}_i(t) = \sum_{j=1}^n W_{ij}(t) \lp \bfx_j(t) + \bfz_{i,j}(t)\rp = \sum_{j=1}^n W_{ij}(t)  \bfx_j(t) + \sum_{j=1}^n W_{ij}(t)  \bfz_{i,j}(t).  
	\end{align*}
	Thus, we get $\nt_i(t) = \sum_{j=1}^n W_{ij}(t) \bfz_{i,j}(t)$ implying 
	\[\EE{\nt_i(t)| \mathcal{F}_t}=0,\quad {\EE{\|\nt_i(t) \|^2 |\mathcal{F}_t} = \zeta^2 \sum_{j=1}^n W^2_{ij}(t) \leq \zeta^2}\]. 
	Hence, the conditions of Assumption~\ref{assum:neighbor} are satisfied. 
\end{example}
  \section{Experimental Results}\label{sec:expr-res}
   In the following, we discuss two sets of simulations to  verify our theoretical results. First, we perform Algorithm~\ref{alg:1} for a \textit{fixed} network to compare fixed step-sizes~\cite{reisizadeh2019robust} versus  diminishing step-sizes (this work). Then, we compare the performance of our algorithm on fixed versus time-varying graph. Throughout this section, we use the unbiased stochastic quantizer in~\eqref{eq:unb-q} with various total number of quantization levels $s$. 
   For CNN experiments over CIFAR-10 and regularized logistic regression problem over MNIST dataset, the mini-batch size is set to be $10$ and $20$, respectively, the loss function is the cross-entropy function, we use $N=10,000$ data points which are distributed across $n=20$ agents. Moreover, for each set of experiments, we distribute the data points across the nodes according to $r_i=p_i/\sum_{i=1}^{20} p_i$, where $p_i$ is drawn uniformly at random from the interval $(0.01,0.9)$. Our codes are implemented using Python and tested on a 2017 MacBook Pro with 16 GB memory and a 2.5 GHz Intel Core i7 processor.
   
    \subsection{\texttt{DIMIX} vs.\ \texttt{QuanTimed-DSGD} over Fixed Network}
    We use regularized logistic regression to validation  our algorithm on the benchmark dataset MNIST.
    \noindent$\triangleright$ \textbf{Data and Experimental Setup.}
   In this  experiments we train a regularized logistic regression to classify the MNIST dataset over a fixed network. First, we generate a random (connected) undirected Erd\"os-Renyi graph with the edge probability $p_c=0.3$, which is fixed throughout this experiment. We also generate a {doubly} stochastic weight matrix ${A:=I-(d_{\max}+1)^{-1}L_G}$ where $L_G$ is the Laplacian matrix  and $d_{\max}$ is the maximum degree of the graph. Finally,  we use the time-invariant weight sequence $W=I+ c\hat{\bfr}_{\min}\text{diag}^{-1}(\bfr)(A-I)$, where ${\hat{\bfr}_{\min} = \min_{i\in[n]} r_i/(1-A_{ii})}$, $c=0.95$, and $\text{diag}(\bfr)$ is a diagonal matrix with $r_i$ as its $i$th diagonal element. It can be verified that $W$ satisfies the conditions of Assumption~\ref{asm:W}, i.e., it is row-stochastic and  satisfies  $\bfr^TW=\bfr^T$. 
    We implemented our algorithm with quantizer parameter $s=3$ and tuned the step-size parameters $(\azr,\nu^\star) = (0.005,1/6)$, $(\bzr,\mu^\star) = (0.6,1/2)$, with $\tu = 2000$. For the fixed step-sizes, we used the termination time $T=7500$, which results in $\alpha(t)=0.001$ and $\beta(t)=0.01$ for all $t\geq 1$.
    \begin{figure}[h]
    \vspace{-3pt}
    \centering
    \includegraphics[width=0.5\linewidth]{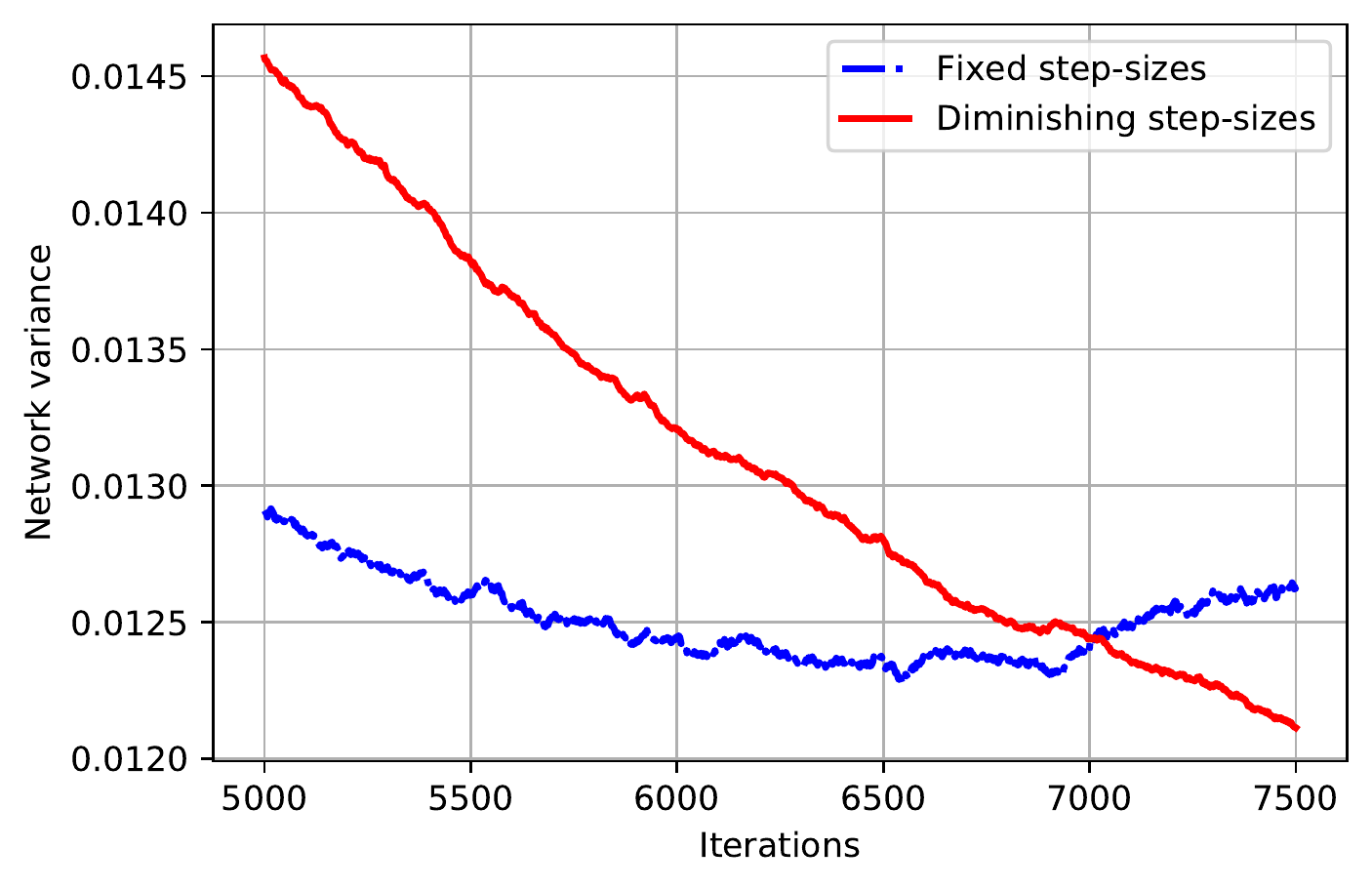}
    \vspace{-3pt}
    \caption{Logistic Regression on MNIST: network variance for fixed and vanishing step-sizes.}\vspace{-3mm}
    \label{fig:mnist-logr-fx-g}
    \end{figure}
    
    \noindent\textbf{Results.} The plot in Figure~\ref{fig:mnist-logr-fx-g} shows \textit{network variance} defined by $\Nr{X(t)-\ones\mx(t)}^2$ for fixed and diminishing step-sizes. It can be observed that with a fixed step-sizes algorithm 
    we can only reach a  \textit{neighborhood} of the average model, but a \textit{consensus} is not necessarily achieved. However, using diminishing step-sizes guarantees that each node's state converges to the average model.
    
    \subsection{Diminishing Step-sizes over Time-varying  vs. Fixed Network}
    Here, we provide simulation results to demonstrate the performance of \texttt{DIMIX} on time-varying networks. For this, we conduct simulations on non-convex (CNN) and convex (linear regression) problems. 
    For the sake of comparison, we apply our algorithm on a corresponding fixed network.  
    
    \noindent\textbf{Data and Experimental Setup.}
    For the non-convex setting, we study the problem of learning the parameters of the LeNet-5 CNN~\cite{lecun1998gradient} for CIFAR-10 dataset. The LeNet-5 CNN architecture consists  of two convolutional layers, two sub-sampling (pooling) layers, and three fully connected layers, with the ReLU activation function. As before we have $n=20$ nodes and $N=10,000$ data points, which are distributed among the agents according to some stochastic vector $\bfr$. 
    For the convex setting, we consider a $100$-dimensional linear regression problem. We generate $N = 300$ data points $\{\xi_1,\ldots, \xi_{300}\}$, where $\xi_i=(\bfu_i,v_i)$, with  $v_i = \bfu_i^T \tilde{\bfx} + \epsilon_i$, and $\bfu_i,\tilde{\bfx}\in \R^{100}$ for all $i\in[300]$.  In order to generate the data, we uniformly and independently draw the entries of each $\bfu_i$ and each  $\epsilon_i$ from  $(0,1)$  and $(0,0.1)$, respectively. Similarly, entries of $\tilde{\bfx}$ are sampled uniformly and independently from $(-1,1)$. The goal is to estimate the unknown coefficient vector $\tilde{\bfx}$, which leads to solving the minimum mean square error problem, i.e.,
    \[f(\bfx):=\min_{\bfx\in \mathbb{R}^d}\frac{1}{2N} \sum_{i=1}^{300} \|v_i - \bfu_i^T \bfx\|^2\].

    \noindent\textbf{Experiments with Fixed Graph.}
     We consider a fixed undirected cyclic graph $\mathcal{G}^C = ([n], \mathcal{E})$, where ${\mathcal{E} = \{(\la i \ra,\la i+1\ra): i\in [n]\}}$, and $\la i \ra := (i-1 \mod n)+1$.  
     For these experiments, the stochastic matrix sequence $W(t)=W$ for any $t$ is given by
     \begin{align}\label{W:fdx}
          W_{ij} = \begin{cases}
          \frac{r_{\la j \ra}}{2(r_{\la i \ra} +r_{\la j \ra})} & j \in\{\la i-1 \ra, \la i+1 \ra\} \\
          \frac{r_{\la i \ra}}{2(r_{\la i \ra}+r_{\la i+1 \ra})} + \frac{r_{\la i \ra}}{2(r_{\la i \ra}+r_{\la i-1 \ra})} & j = i \\
          0 & \textrm{otherwise}.
          \end{cases}
     \end{align}
     
    \noindent\textbf{Experiments with Time-varying Graph.}
    To evaluate Algorithm~\ref{alg:1} for a time-varying network, we use cyclic gossip  \cite{boyd2006randomized,dimakis2010gossip,kar2008distributed,ram2010asynchronous,lee2015asynchronous} iterations. Here, the algorithm goes through the cycle graph $\mathcal{G}^C$ (described above), and a pair of  neighbors (on the cycle) exchange  their information at a time. More precisely, at time $t$, the averaging graph $\mathcal{G}(t) = ([n], \mathcal{E}(t))$ only  consists of a single edge $\mathcal{E}(t) = \{({\la t \ra}, {\la t+1 \ra} )\}$. For $\clg(t)$, we let 
    \begin{align}\label{W:var}
        [W(t)]_{ij} =
        \begin{cases}
        \frac{r_{\la j\ra}}{ r_{{\la t\ra}} + r_{{\la t+1\ra}}} & i,j\in \{{\la t\ra}, {\la t+1\ra}\}\\
        1 & i=j \notin \{{\la t\ra}, {\la t+1\ra}\}\\
        0 & \textrm{otherwise}.
        \end{cases}
    \end{align}
    Note that each edge in $\clg^C$ will be visited periodically, and hence, the sequence of stochastic matrices $\{W(t)\}$  in~\eqref{W:var} is $B$-connected with $B=n$. 
    \begin{figure}[h]
    \centering
    \includegraphics[width=0.49\linewidth]{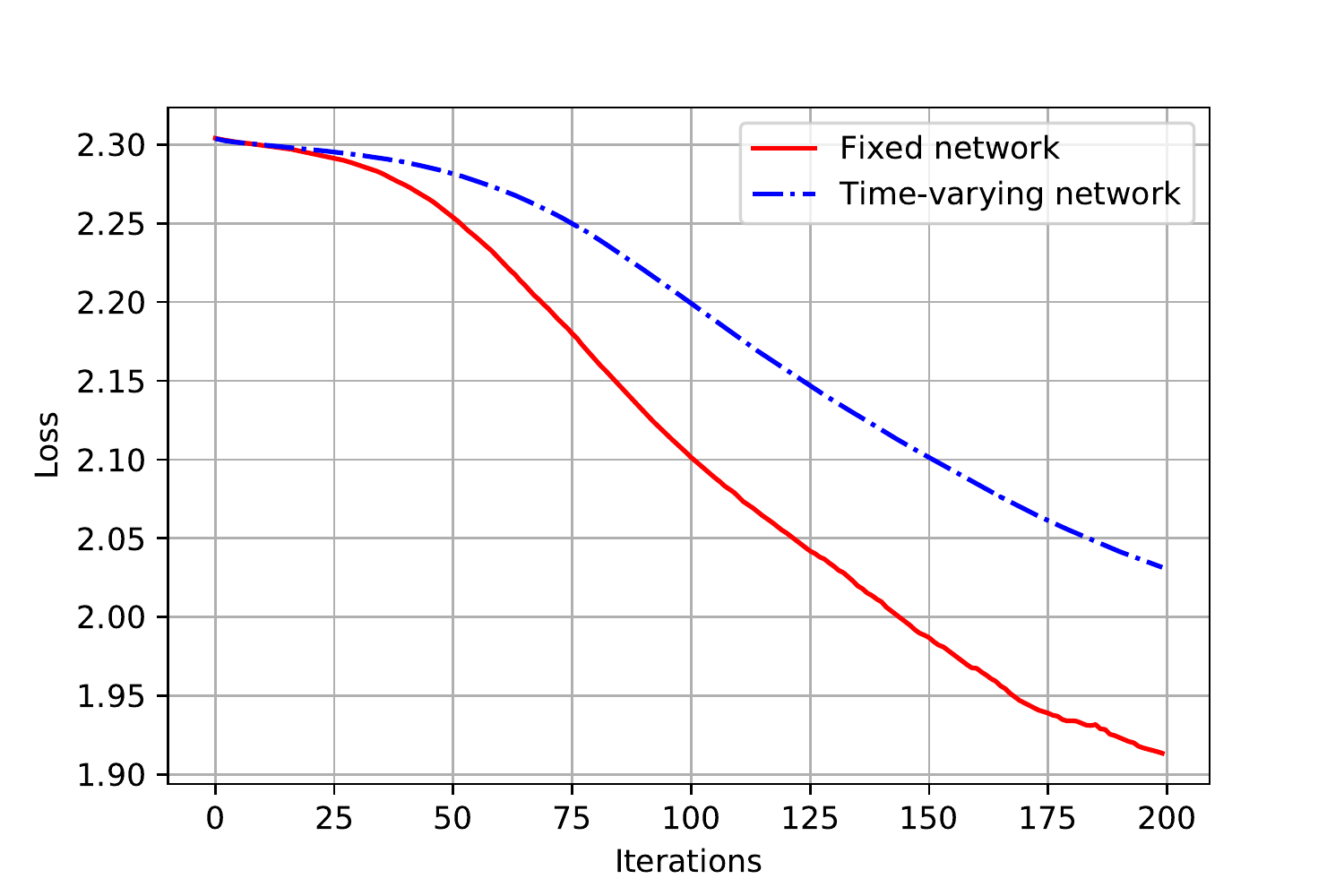}
    \hfill
    \includegraphics[width=0.5\linewidth]{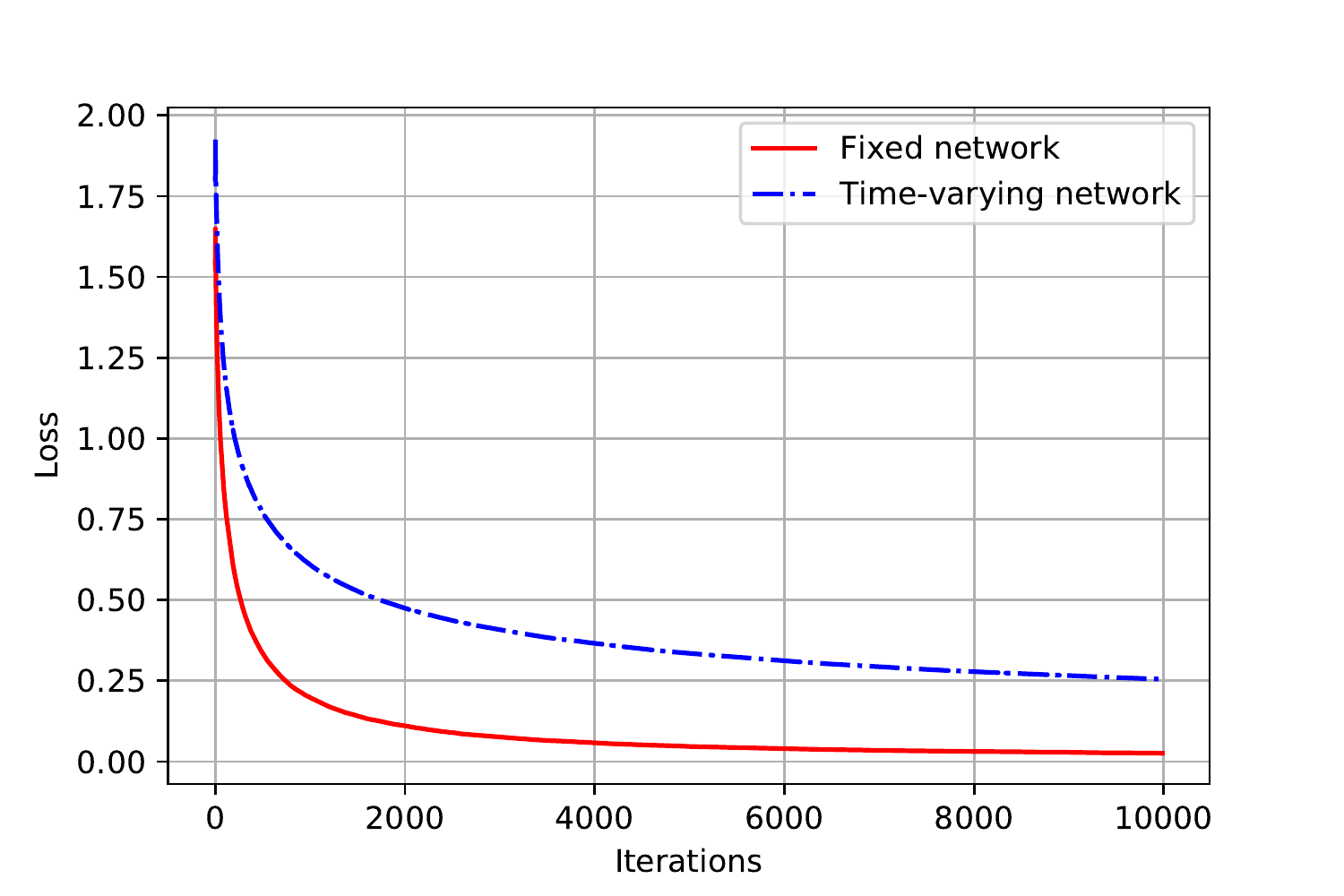}
    \caption{Training Loss vs. Iterations: LeNet-5 on CIFAR-10 (left), and Linear Regression (right).}
    \label{fig:time-varying-g}
    \end{figure}
    
     \noindent\textbf{Results.}
     In Figure~\ref{fig:time-varying-g}, the plot on the left demonstrates the training loss of LeNet-5 for {CIFAR-10} dataset and the plot on the right shows the training loss of the linear regression. 
     Here, `Fixed network' refers to the full cycle with stochastic matrix  in~\eqref{W:fdx} and `Time-varying network' refers to the network with stochastic matrix sequence in~\eqref{W:var}. For the CNN, the parameters of the dynamics in~\eqref{eq:upd_nd} are fine-tuned to ${(\azr,\nu^\star) = (0.02,1/6)}$, ${(\bzr,\mu^\star) = (0.05,1/2)}$,  $\tu=0$, and  $s=6$ is the parameter of the stochastic quantizer for both  networks.  For the linear regression model, the parameter are adjusted to $(\azr,\nu) = (6,1/4)$, $(\bzr,\mu) = (16,3/4)$, and $\tu=1500$ with the quantizer parameter $s=6$.
     It can be observed that the algorithm works for both fixed and time-varying networks. However, the performance of the algorithm over the time-varying network naturally suffers from a slower convergence, which is due to a slower mixing of information over the network.
     
    \section{Auxiliary Lemmas}\label{sec:aux-lemmas}
    The following lemmas are used to facilitate the proof of the main result of the paper. We present the proofs of Lemmas~\ref{lm:nor_squr}--\ref{lemma:transition2} in Section~\ref{sec:proof-aux}. 
\begin{lemma}\label{lm:nor_squr}
   For any pair of vectors $\bfu$, $\bfv$, and any scalar $\omega>0$, we have 
   \begin{align}
        \|\bfu+\bfv\|^2 &\leq (1+\omega)  \|\bfu\|^2 + \lp 1+ \frac{1}{\omega} \rp \|\bfv\|^2.
   \end{align}
   Similarly, for matrices $U$ and $V$ and any scalar $\omega>0$, their $\bfr$-norm satisfy
   \begin{align}
       \|U+V\|_\bfr^2 \leq (1+\omega) \|U\|^2_\bfr + \lp 1+\frac{1}{\omega}\rp  \|V\|^2_\bfr.
   \end{align}
\end{lemma}
\begin{lemma}\label{lemma:Psum}
Let $\{\beta(t)\}$ be a sequence in $\R$ and $\lambda$ be a positive scalar. Then the following identities hold
\begin{align}\label{eq:esumprod}
\sum_{s=1}^{t-1} \beta(s)\!\!\prod_{k=s+1}^{t-1}(1-\lambda \beta(k))&=\frac{1}{\lambda} -\frac{1}{\lambda}\prod_{k=1}^{t-1}(1-\lambda \beta(k)),\quad \mbox{for all $1\leq t$},\\ 
    	 \sum_{t=s+1}^{T} \beta(t)\! \prod_{k=s+1}^{t-1}(1-\lambda \beta(k))&=\frac{1}{\lambda}-\frac{1}{\lambda} \prod_{k=s+1}^{T} (1-\lambda \beta(k))\quad \mbox{for all $1\leq s\leq T$}.\label{eq:esumprodb}
    	\end{align}
    	As a result, for any sequence $\{\beta(t)\}$ in $[0,1]$ and $\lambda>0$, we have 
    	\[\sum_{s=1}^{t-1} \beta(s)\prod_{k=s+1}^{t-1}(1-\lambda \beta(k))\leq  \frac{1}{\lambda},\]
    	and 
    	\[\sum_{t=s+1}^{T} \beta(t)\prod_{k=s+1}^{t-1}(1-\lambda \beta(k))\leq \frac{1}{\lambda}.\]
    \end{lemma}
    
    \begin{lemma}\label{lm:sum-exp}
    For any $\dl\in \mathbb{R}$, $\tau \geq 0$, and $T\geq 1$, we have
    \begin{align}
    \sum_{t=1}^T (t+\tau)^{\dl} \leq
    \begin{cases}
    \frac{\tau^{1+\delta}}{|1+\dl|} & \textrm{if $\dl<-1$},\\
    \ln\lp \frac{T}{\tau} +1 \rp& \textrm{if $\dl=-1$},\\
    \frac{2^{1+\dl}}{1+\dl} (T+\tau)^{1+\dl} & \textrm{if $\dl>-1 $}.
    \end{cases}
    \end{align}
    \end{lemma}

	\begin{lemma}\label{lemma:matrixprop}
	    For an $n\times m$ matrix $A$ and $m\times q$ matrix $B$, we have:
	    \begin{align*}
             \lnr AB \rnr_{\bfr}	   \leq     \lnr A \rnr_{\bfr} \lnr B \rnr_{F}.  
	    \end{align*}
	\end{lemma}

As we discussed in Remark~\ref{rem:vanishing-graph}, we cannot use the conventional results (e.g., in~\cite{demarzo2003persuasion,nedic2008distributed,nedic2009distributed,nedic2013distributed,touri2013product,touri2014endogenous,yuan2016convergence,tatarenko2017non,aghajan2020distributed}) to bound the norm of a vector after averaging with matrices with \textit{diminishing} weights. The following lemma provides a new bounding techniques, which we will use in the proof of the main result of this work.

\begin{lemma}\label{lemma:transition2}
	    Let $\Wt$ satisfy the Connectivity Assumption \ref{asm:W} with parameters $(B,\minW)$,  and let  $\{A(t)\}$ be given by $A(t)=(1-\beta(t))I+\beta(t)W(t)$ where $\beta(t)\in (0,1]$ for all $t$ and $\{\beta(t)\}$ is a non-increasing  sequence. Then, for any matrix $U\in \mathbb{R}^{n\times d}$, parameter $\lambda=\frac{\minW \rmin}{2Bn^2}$, and all $t>s\geq 1$, we have 
	    \[\Nr{\lp A(t-1)A(t-2)\cdots A(s+1)-\ones \bfr^T\rp U }^2 \leq \kap \prod_{k=s+1}^{t-1}(1-\lambda\beta(k))\Nr{U}^2,\]
	    where $\kap = \lp 1-B\lambda \bzr\rp^{-1}$ and $\bzr=\beta(1)$. 
	\end{lemma}

    \section{Proof of Theorem~\ref{thm:non_cvx}}\label{sec:proof_nn_cvx}
    For our analysis, first we obtain an expression for the average reduction of the objective function $f(\cdot)$ at the average of the states, i.e., $\mx(t) = \bfr^T X(t) =\sum_{i=1}^n r_i \bfx_i(t)$. Recall that $\bfr^T W(t) = \bfr^T$ for all $t\geq 1$. Hence, multiplying both sides of~\eqref{eq:dyn-col-mat} by $\bfr^T$, we get
    \begin{align*}
     \mx(t+1) = \mx(t)+\et\bfr^T \ER{t} - \at\et\bfr^T\nabla f(X(t)).  
    \end{align*}
    From Assumption~\ref{asm:f}-(a)
    Lemma~3.4 in \cite{bubeck2014convex} we can conclude
    \begin{align*}
        f(\mx(t+1)) - f(\mx(t)) - \la \nabla  f(\mx(t)), \mx(t+1) -\mx(t)\ra  \leq \frac{\li}{2} \lnr \mx(t+1) -\mx(t) \rnr^2,
    \end{align*}
    or equivalently, 
    \begin{align*}
        f(\mx(t+1)) \leq &f(\mx(t)) + \beta(t) \la \nabla  f(\mx(t)), \bfr^T E(t) \ra 
        -\alpha(t) \beta(t) \la \nabla  f(\mx(t)), \bfr^T \nabla f(X(t))\ra \cr
        &+ \beta(t)^2 \frac{\li}{2} \lnr \bfr^T E(t) - \alpha(t) \bfr^T \nabla f(X(t)) \rnr^2.
    \end{align*}
    Then, since $\{X(t)\}$ is adapted to the  filtration $\{\F_t\}$ and using Assumption~\ref{assum:neighbor}, we arrive at
    \begin{align}\label{eq:Colbnd}
        \EE{f\left(\mx(t+1)\right)\middle|\F_t}&\leq f\left(\mx(t)\right) - \at\et\la\nabla f\left(\mx(t)\right),\bfr^T\nabla f(X(t))\ra \cr
        &\phantom{--} + \beta^2(t)\frac{\li}{2}\EE{ \left\|\bfr^T\left[\ER{t}
        - \at\nabla f(X(t))\right]\right\|^2\middle|\F_t}.
    \end{align}
    	Let us focus on the second term in \eqref{eq:Colbnd}. Using the identity $2 \la \bfx - \bfy \ra  = \lnr \bfx \rnr^2 + \lnr \bfy \rnr^2 - \lnr \bfx - \bfy \rnr^2$, we can write
	\begin{align}\label{eq:Colbnd-2-1}
	    \la\nabla f\left(\mx(t)\right),\bfr^T\nabla f(X(t))\ra \!=\! \| \nabla f\lp \mx(t)\rp\|^2 \!+\! \|\bfr^T\nabla f(X(t))\|^2\!-\! 2  \| \nabla f\lp\mx(t)\rp - \bfr^T\nabla f(X(t))\|^2.
	\end{align}
	Moreover, we have 
	\begin{align*}
	     \| \nabla f\lp\mx(t)\rp - \bfr^T\nabla f(X(t))\|^2
	    &= \lnr \sum_{i=1}^n r_i\nabla f_i(\mx(t)) - \sum_{i=1}^n r_i \nabla f_i(\bfx_i(t)) \rnr^2\cr 
	    &= \lnr \sum_{i=1}^n r_i\left(\nabla f_i(\mx(t)) - \nabla f_i(\bfx_i(t))\right) \rnr^2\cr 
	    &\leq \sum_{i=1}^n r_i\lnr\nabla f_i(\mx(t)) - \nabla f_i(\bfx_i(t)) \rnr^2,
	\end{align*}
	where the last inequality holds as $\|\cdot\|^2$ is a convex function and the $\bfr$ is a stochastic vector. 
	Therefore, using the above inequality and  Assumption~\ref{asm:f}-(a), we get
	\begin{align}\label{eq:Colbnd-2-2}
	     \| \nabla f\lp\mx(t)\rp - \bfr^T\nabla f(X(t))\|^2
	    &\leq \sum_{i=1}^n r_i\lnr\nabla f_i(\mx(t)) - \nabla f_i(\bfx_i(t)) \rnr^2\cr 
	    &\leq \li^2\sum_{i=1}^n r_i\lnr \mx(t) - \bfx_i(t) \rnr^2.
	\end{align}
	
    Next, we analyze the last term in~\eqref{eq:Colbnd}. 
    From Assumption~\ref{assum:neighbor} we have  $\EE{E(t)|\F_t}=0$, which leads to  
    \begin{align}\label{eq:Colbnd2}
        \EE {\left\|\bfr^T\left[\ER{t}
        - \at\nabla f(X(t))\right]\right\|^2\middle|\F_t} =\EE{\left\|\bfr^T\ER{t} \right\|^2\middle|\F_t}+ \left\|
         \at\bfr^T\nabla f(X(t))\right\|^2.
    \end{align}
    For the first term in~\eqref{eq:Colbnd2}, we again exploit Assumption~\ref{assum:neighbor} which implies
  \begin{align}
        \EE{\left\|\bfr^T \ER{t}
        \right\|^2 \middle|\F_t}
        & = \bfr^T\EE{\ER{t}\ER{t}^T
         \middle|\F_t}  \bfr
        \nonumber\\
        &\leq \bfr^T (\sm \ones \ones^T) \bfr = \gamma,
        \label{eq:quant-noise}
   \end{align}
   where we used the fact that $\bfr^T\ones =1$ and the inequality holds since
   \begin{align*}
       \lb\EE{\ER{t}\ER{t}^T \middle|\F_t}\rb_{ij}&=\EE{\nt_i(t) \nt_j(t) | \F_t} \cr 
        &\leq \sqrt{\EE{\|\nt_i(t)\|^2  | \F_t} \EE{\| \nt_j(t)\|^2 | \F_t} }\leq \gamma,
   \end{align*}
   for all $i,j \in\  [n]$. Thus, the last term in~\eqref{eq:Colbnd} is upper bounded as 
   \begin{align}\label{eq:Colbnd-3}
       \EE{ \left\|\bfr^T\left[\ER{t}
        - \at\nabla f(X(t))\right]\right\|^2\middle|\F_t} \leq \gamma + \alpha^2(t) \| \bfr^T\nabla f(X(t))\|^2.
   \end{align}
	Therefore, replacing \eqref{eq:Colbnd-2-1}, \eqref{eq:Colbnd-2-2}, and \eqref{eq:Colbnd-3} in \eqref{eq:Colbnd} we get
	\begin{align*}
        &\EE{f\left(\mx(t+1)\right)\middle|\F_t}\cr
        &\leq f\left(\mx(t)\right) - \at\et\la\nabla f\left(\mx(t)\right),\bfr^T\nabla f(X(t))\ra \cr
        &\quad + \beta^2(t)\frac{\li}{2}\EE{ \left\|\bfr^T\left[\ER{t}
        - \at\nabla f(X(t))\right]\right\|^2\middle|\F_t}\cr 
        &\leq f\lp \mx(t)\rp - \frac{1}{2}\at\et\lp\lnr \nabla f(\mx(t))\rnr^2 +\lnr \bfr^T\nabla f(X(t))\rnr^2-2 \li^2\sum_{i=1}^n r_i\lnr \mx(t)-\bfx_i(t)\rnr^2 \rp\cr
        &\quad + \beta^2(t)\frac{\li}{2}\lp\gamma+\alpha^2(t)\lnr\bfr^T\nabla f(X(t))\rnr^2\rp\cr
        &= f\lp \mx(t)\rp -\frac{1}{2}\at\et\lnr\nabla f(\mx(t))\rnr^2 -\frac{1}{2} \at\et\lp 1-\at\et \li\rp\lnr \bfr^T \nabla f(X(t))\rnr^2  \cr 
        &\quad + \at\et \li^2\Nr{X(t)  - \ones\mx(t)}^2 + \beta^2(t)\frac{\li}{2}\gamma.
    \end{align*}
	Taking the expectation of both sides lead to
	\begin{align}\label{eq:main-ineq}
        \EE{f\left(\mx(t+1)\right)}
        &\leq \EE{ f\lp \mx(t)\rp} -\frac{1}{2}\at\et\EE{\lnr\nabla f(\mx(t))\rnr^2}\cr
        &\quad -\frac{1}{2} \at\et\lp 1-\at\et \li\rp\EE{\lnr \bfr^T \nabla f(X(t))\rnr^2}  \cr 
        &\quad + \at\et \li^2\EE{\Nr{X(t) - \ones \mx(t)}^2} + \beta^2(t)\frac{\li}{2}\gamma.
    \end{align}
\subsection{State Deviations from the Average State: \texorpdfstring{$\mathbb{E} \left[ \left\| X(t) - \ones\bar{\bfx}(t)\right\|_{\bfr}^2\right]$}{}}
	\label{sec:node-to-avg}
	
	Note that the dynamics in~\eqref{eq:dyn-col-mat} can be viewed as the linear time-varying system  
	\begin{align}\label{eq:controlled}
	    X(t+1)=A(t)X(t)+U(t),
	\end{align}
    with
    \begin{align*}
    A(t)&=((1-\et)I+\et W(t)),\cr
    U(t)&=\et\ER{t} - \at\et\nabla f(X(t)). 
    \end{align*}
     The solution of~\eqref{eq:controlled} is given by
	\begin{align}
	    X(t)=\sum_{s=1}^{t-1}\Phi(t:s)U(s) + \Phi(t:0)X(1),
	    \label{eq:Z(t)}
	\end{align}
	where $\Phi(t:s)=A(t-1)\cdots A(s+1)$ with $\Phi(t:t-1)=I$ is the transition matrix of the linear system \eqref{eq:controlled}. 
	For the simplicity of notations, let us define 
	\begin{align*}
	    P(t:s)&:=\beta(s)(\Phi(t:s)-\ones\bfr^T) = \beta(s)\lp A(t-1)\cdots A(s+1)-\ones\bfr^T\rp . 
	\end{align*}
     As a result of Lemma~\ref{lemma:transition2}, we have
	     \[\Nr{ P(t:s)U}\leq \pi(t:s)\Nr{U},\]
	     where $\pi(t:s)$ is defined by 
	     \begin{align}\label{eqn:pidef}
	         \pi(t:s):=\beta(s)\kap^{\frac{1}{2}}\prod_{k=s+1}^{t-1}(1-\lambda\beta(k))^{\frac{1}{2}}.
	     \end{align}
	Now consider the dynamic in~\eqref{eq:Z(t)}. Assuming $X(1)={\mathbf{0}}$, we get
	\begin{align}
	    X(t)=\sum_{s=1}^{t-1}\Phi(t:s)U(s).
	    \label{eq:Z2(t)}
	\end{align}
	Moreover,  multiplying both sides of \eqref{eq:Z2(t)} from the left by $\bfr^T$ and using $\bfr^TA(t)=\bfr^T$, we get 
	\begin{align}
	\mx(t)=\bfr^T X(t) = \sum_{s=1}^{t-1}\bfr^T \Phi(t:s)U(s) = \sum_{s=1}^{t-1}\bfr^T U(s).
	\label{eq:Z:mean}
	\end{align}
	Then, subtracting~\eqref{eq:Z:mean} from~\eqref{eq:Z2(t)}, and plugging the definition of $U(s)$ we have
	\begin{align*}
	    X(t)-\ones \bar{\bfx}(t)
	    &= \sum_{s=1}^{t-1}\Phi(t:s)U(s) - \sum_{s=1}^{t-1} \ones \bfr^T U(s)\cr
	    &=\sum_{s=1}^{t-1}(\Phi(t:s)-\ones\bfr^T)U(s)\cr
	    &=\sum_{s=1}^{t-1}\beta(s)(\Phi(t:s)-\ones\bfr^T)\bigg[ \ER{s} - \alpha(s)\nabla f(X(s))\bigg]\cr
	    &=\sum_{s=1}^{t-1} P(t:s)\ER{s} -\sum_{s=1}^{t-1}\alpha(s)P(t:s)\nabla f(X(s)).
	\end{align*}
	Therefore, using 
	Lemma~\ref{lm:nor_squr} with $\omega=1$, we get
    \begin{align*}
	    \Nr{X(t)-\ones\bar{\bfx}(t)}^2	    &=\Nr{\sum_{s=1}^{t-1} P(t:s)\ER{s}-\sum_{s=1}^{t-1}\alpha(s)  P(t:s)\nabla f(X(s))}^2 \cr
	    & \leq 2\Nr{\sum_{s=1}^{t-1} P(t:s)\ER{s}}^2 +2\Nr{\sum_{s=1}^{t-1}\alpha(s)  P(t:s)\nabla f(X(s))}^2.
	\end{align*}
	By expanding $\Nr{\sum_{s=1}^{t-1}  P(t:s)\ER{s}}^2$, we get
	\begin{align}\label{eq:Mtmain}
	     \Nr{X(t)-\ones\bar{\bfx}(t)}^2&=2\sum_{s=1}^{t-1}  \Nr{ P(t:s)\ER{s}}^2 + 2\sum_{s\neq q} \la P(t:s)\ER{s}, P(t:q)\ER{q}\ra \cr
	    &\phantom{=} + 2\Nr{ \sum_{s=1}^{t-1} \alpha(s)   P(t:s)\nabla f(X(s))}^2.
	\end{align}
	Recall from  Assumption~\ref{assum:neighbor} that $\EE{E(q)|\F_q}=0$. Moreover, since  $E(s)$ is measurable with respect to $\F_q$ for  $q>s$, we have
	\begin{align*}
        \EE{\la P(t:s) \ER{s}, P(t:q) \ER{q}\ra} & =	
        \EE{\EE{\la P(t:s) \ER{s}, P(t:q) \ER{q}\ra}\md \F_q}\cr 
        & =\EE{\la P(t:s) \ER{s}, P(t:q) \EE{\ER{q}\md \F_q}\ra}=0. 
    \end{align*}
    Using a similar argument for $q<s$ and conditioning on $\F_s$, we conclude that the above relation holds for all $q\not=s$. Therefore, taking the expectation of both sides of \eqref{eq:Mtmain} and noting that the average of the second and the forth terms are zero,  we get
	\begin{align}
	    \EE{\Nr{X(t)-\ones \bar{\bfx}(t)}^2}
	    &\leq 2\sum_{s=1}^{t-1} \EE{ \Nr{ P(t:s)\ER{s}}^2}  +2\EE{ \Nr{\sum_{s=1}^{t-1} \alpha(s)  P(t:s)\nabla f(X(s))}^2} 
	    \label{eq:Mtmain:exp} 
	\end{align}
	 We continue with bounding  the first term in~\eqref{eq:Mtmain:exp}. First note that Assumption~\ref{assum:neighbor} implies
	 \begin{align}
	     \EE{\Nr{E(s)}^2} = \EE{\EE{\Nr{E(s)}^2|\F_s}} = \EE{ \sum_{i=1}^n r_i \EE{\|\bfe_i(s)\|^2 |\F_s}} \leq \EE{\sum_{i=1}^n r_i \gamma} = \gamma.
	 \end{align}
	 This together with Lemma~\ref{lemma:transition2} leads to 
   \begin{align}\label{eq:X-mx-term}
   \sum_{s=1}^{t-1} \EE{ \Nr{ P(t:s)\ER{s} }^2}
   & \leq \lb \sum_{s=1}^{t-1} \beta^2(s)\kap \prod_{k=s+1}^{t-1} (1-\lambda \beta(k)) \EE{ \Nr{ \ER{s}}^2} \rb \cr   
   & \leq \sm \kap \sum_{s=1}^{t-1} \lb \beta^2(s) \prod_{k=s+1}^{t-1} (1-\lambda \beta(k)) \rb.
   \end{align}
    To bound the second term in~\eqref{eq:Mtmain:exp}, we use the triangle inequality for the norm  $\Nr{\cdot}$, and write 
\begin{align}
	\begin{split}\label{eq:Egradient}
	&\EE{ \Nr{\sum_{s=1}^{t-1} \alpha(s)  P(t:s)\nabla f(X(s))}^2}\cr
	&\leq \EE{ \lp\sum_{s=1}^{t-1} \Nr{\alpha(s)  P(t:s)\nabla f(X(s))}\rp^2}\cr 
	&=\sum_{1\leq s,q \leq t-1} \EE{ \alpha(s) \Nr{ P(t:s)\nabla f(X(s))} \alpha(q) \Nr{ P(t:q)\nabla f(X(q))}}.
	  \end{split}
	  \end{align}
	  Using Lemma~\ref{lemma:transition2} and $2ab\leq a^2+b^2$, we can upper-bound this expression as
	 \begin{align}
	 \sum_{1\leq s, q \leq t-1} & \EE{ \alpha(s) \Nr{ P(t:s)\nabla f(X(s))} \alpha(q) \Nr{ P(t:q)\nabla f(X(q))}}\cr
	  & \leq \sum_{1\leq s, q\leq t-1} 
	  \EE{\alpha(s) \pi(t:s) \Nr{\nabla f(X(s))} \alpha(q) \pi(t:q) \Nr{\nabla f(X(q))}}  \cr
	  &= \sum_{1\leq s, q\leq t-1}  \pi(t:s) \pi(t:q) \EE{ \alpha(s)\Nr{\nabla f(X(s))}  \alpha(q) \Nr{\nabla f(X(q))}} \cr 
	  &\leq \frac{1}{2} \sum_{1\leq s, q\leq t-1}  \pi(t:s) \pi(t:q) \EE{ \alpha^2(s)\Nr{\nabla f(X(s))}^2+ \alpha^2(q) \Nr{\nabla f(X(q))}^2} \cr
	  &= \sum_{1\leq s, q\leq t-1}  \pi(t:s) \pi(t:q) \EE{\alpha^2(s) \Nr{\nabla f(X(s))}^2}\cr 
	  &= \lp \sum_{q=1}^{t-1} \pi(t:q)\rp 
	\lp \sum_{s=1}^{t-1}  \alpha^2(s) \pi(t:s) \EE{\Nr{\nabla f(X(s))}^2}\rp,
	  \label{eq:term4-2}
	\end{align}
	where $\pi(t:s)$ is given by \eqref{eqn:pidef}. Next, the fact that $\sqrt{1-x} \leq 1-x/2$ and Lemma~\ref{lemma:Psum} imply  
	\begin{align}
    	\sum_{q=1}^{t-1} \pi(t:q) &= \sum_{q=1}^{t-1}  \lb \beta(q)\kap^{\frac{1}{2}}   \prod_{k=q+1}^{t-1}  (1-\lambda \beta(k))^{\frac{1}{2}} \rb\cr 
    	& \leq \sum_{q=1}^{t-1} \beta(q)\kap^{\frac{1}{2}}   \prod_{k=q+1}^{t-1} \lp 1-\frac{\lambda}{2} \beta(k)\rp \cr
    	&\leq \frac{2}{\lambda}\kap^{\frac{1}{2}}.
    	 \label{eq:sum-Aq}
	\end{align}
	Using this inequality in~\eqref{eq:term4-2}, we get
    \begin{align}\label{eq:Egradient-cross}
	\EE{ \Nr{\sum_{s=1}^{t-1} \alpha(s)  P(t:s)\nabla f(X(s))}^2} 
	&\leq  \frac{2}{\lambda}\kap^{\frac{1}{2}} \sum_{s=1}^{t-1} \lb \alpha^2(s) \pi(t:s)   \EE{\Nr{\nabla f(X(s))}^2}\rb.
	\end{align}
	Finally, using the bounds obtained  in~\eqref{eq:X-mx-term}  and~\eqref{eq:Egradient-cross} in~\eqref{eq:Mtmain:exp} we arrive at 
	\begin{align}\label{eq:var}
	   \EE{\Nr{X(t)-\ones\bar{\bfx}(t)}^2}
	    &\leq 
	     2\sm\kap \sum_{s=1}^{t-1} \lb \beta^2(s) \prod_{k=s+1}^{t-1} (1-\lambda \beta(k)) \rb
	     \cr
	     &\phantom{\leq} +\frac{4}{\lambda}\kap^{\frac{1}{2}} \sum_{s=1}^{t-1} \lb \alpha^2(s) \pi(t:s) \EE{\Nr{\nabla f(X(s))}^2}\rb.
	\end{align}

\subsection{Analysis of the overall deviation: \texorpdfstring{$\sum_{t=1}^T \alpha(t)\beta(t) \EE{\Nr{X(t)-\ones\bar{\bfx}(t)}^2}$}{A22}}
Our goal here is to bound the overall weighted deviation of the states from their average. 
First recall the bound for $\EE{\Nr{X(t)-\ones\bar{\bfx}(t)}^2}$, derived in Section~\ref{sec:node-to-avg} for each $t$. Our goal here is to bound 
\begin{align}
    \sum_{t=1}^T \alpha(t)\beta(t) \EE{\Nr{X(t)-\ones\bar{\bfx}(t)}^2}
    \leq 
	     &2\sm\kap \sum_{t=1}^T \alpha(t)\beta(t) \sum_{s=1}^{t-1} \lb \beta^2(s)\!\! \prod_{k=s+1}^{t-1}\! (1-\lambda \beta(k)) \rb
	     \label{eq:var-sum}\\
	     & +\frac{4}{\lambda}\kap^{\frac{1}{2}} \sum_{t=1}^T \alpha(t)\beta(t) \sum_{s=1}^{t-1} \lb \alpha^2(s) \pi(t:s) \EE{\Nr{\nabla f(X(s))}^2}\rb.\nonumber
\end{align}
Focusing on the first term in~\eqref{eq:var-sum}, we can write
\begin{align}\label{eq:sum-weight-ave-dist-to mean:1}
 \sum_{t=1}^T\! \lb\alpha(t)\beta(t) \sum_{s=1}^{t-1}\lb \beta^2(s)\!\! \prod_{k=s+1}^{t-1} (1-\lambda \beta(k)) \rb \rb \! & = \!\sum_{s=1}^{T-1} \lb \beta^2(s)\!\! \sum_{t=s+1}^{T} \lb \alpha(t)\beta(t)\! \prod_{k=s+1}^{t-1}  (1-\lambda \beta(k)) \rb  \rb \cr
 &\leq \!\sum_{s=1}^{T-1} \lb  \alpha(s) \beta^2(s)\!\! \sum_{t=s+1}^{T} \lb \beta(t)\! \prod_{k=s+1}^{t-1}  (1-\lambda \beta(k))\rb \rb \cr
 &\leq \frac{1}{\lambda}\sum_{s=1}^{T-1} \alpha(s) \beta^2(s),
\end{align}
where the first inequality is due to the fact that $\alpha(t)\leq \alpha(s)$ for $t>s$, and the second one follows from Lemma~\ref{lemma:Psum}. Similarly, using the fact that $\alpha(t)\leq \alpha(s)$ for $t>s$, for the second term in~\eqref{eq:var-sum}, we have 
\begin{align}
  &\sum_{t=1}^T \lb \alpha(t) \beta(t)  
	       \sum_{s=1}^{t-1} \lb \alpha^2(s) \pi(t:s) \EE{\Nr{\nabla f(X(s))}^2}\rb\rb \cr
    &=  
    \sum_{s=1}^{T-1} \lb \alpha^2(s) \EE{\Nr{\nabla f(X(s))}^2} \sum_{t=s+1}^T 
     \alpha(t) \beta(t) \pi(t:s)
     \rb \cr
    &\leq 
    \sum_{s=1}^{T-1} \lb \alpha^3(s)  \EE{\Nr{\nabla f(X(s))}^2} \sum_{t=s+1}^T 
      \beta(t) \pi(t:s)
     \rb.
\label{eq:sum-weight-ave-dist-to mean:2}
\end{align}
Recall that $\pi(t:s) = \beta(s)\kap^{\frac{1}{2}} \prod_{k=s+1}^{t-1} (1-\lambda \beta(k))^{\frac{1}{2}}$. 
Then, using the fact that $\sqrt{1-x} \leq 1-x/2$, we can write 
\begin{align}
    \sum_{t=s+1}^T 
     \beta(t) \pi(t:s) &= \sum_{t=s+1}^T \lb \beta (t) \beta(s)\kap^{\frac{1}{2}} \prod_{k=s+1}^{t-1} (1-\lambda \beta(k))^{\frac{1}{2}} \rb\cr
    &\leq\beta(s)\kap^{\frac{1}{2}} \sum_{t=s+1}^T \lb \beta (t)  \prod_{k=s+1}^{t-1} \lp 1-\frac{\lambda}{2} \beta(k)\rp \rb 
    \leq \frac{2}{\lambda} \beta(s)\kap^{\frac{1}{2}}, 
    \label{eq:sum-B}
\end{align}
where the last inequality follows from Lemma~\ref{lemma:Psum}. Therefore, from~\eqref{eq:sum-weight-ave-dist-to mean:2} and ~\eqref{eq:sum-B} we have 
\begin{align}
     \sum_{t=1}^T\! \lb\! \alpha(t) \beta(t) \!
	       \sum_{s=1}^{t-1} \! \lb \! \alpha^2(s) \pi(t:s)   \EE{\Nr{\nabla f(X(s))}^2}\rb\rb 
    &\! \leq \! 
    \frac{2}{\lambda}\kap^{\frac{1}{2}}
    \!\sum_{s=1}^{T-1} \! \lb \alpha^3(s)  \beta(s) \EE{\Nr{\nabla f(X(s))}^2}  
     \rb.
     \label{eq:sum-weight-ave-dist-to mean:3}
\end{align}
Therefore, from~\eqref{eq:sum-weight-ave-dist-to mean:1} and~\eqref{eq:sum-weight-ave-dist-to mean:3} we can conclude 
	\begin{align}
	   \sum_{t=1}^T \! \alpha(t) \beta(t) &\EE{\Nr{X(t)-\ones\bar{\bfx}(t)}^2}
	     \! \leq \!
	    \frac{2\sm\kap}{\lambda}\! \sum_{t=1}^{T-1} \! \alpha(t) \beta^2(t) 
	    \!+\! 
	    \frac{8\kap}{\lambda^2}  \sum_{t=1}^{T-1}\! \lb \alpha^3(t) \beta(t)  \EE{ \Nr{\nabla f(X(t))}^2}\rb\nonumber\\
	    & \! \leq \!
	    \frac{2\sm\kap}{\lambda}\! \sum_{t=1}^{T} \! \alpha(t) \beta^2(t) 
	    \!+\! 
	    \frac{8\kap}{\lambda^2}  \sum_{t=1}^{T}\! \lb \alpha^3(t) \beta(t)  \EE{ \Nr{\nabla f(X(t))}^2}\rb.
	    \label{eq:total-dev}
	\end{align} 
\subsection{Bounding \texorpdfstring{$\EE{ \Nr{\nabla f(X(t))}^2}$}{A23}}
    In this part, we study  $\EE{\lnr \nabla f(X(t)\rnr_\bfr}^2$, to provide an upper bound for it.  Following \cite{reisizadeh2019robust}, we can  rewrite $\nabla f(X(t))$ as
	\[\nabla f(X(t))=3 \lb \frac{1}{3}\lp\nabla f(X(t))-\nabla f(\ones \mx(t)\rp+ 
	\frac{1}{3}\lp\nabla f(\ones \mx(t))-\ones\nabla f(\mx(t))\rp
	+\frac{1}{3} \ones\nabla f(\mx(t))\rb.\]  
	Then, since $\lnr \cdot \rnr^2_\bfr$ is a convex function, we have
	\begin{align}\label{eq:norm-grad-f-X}
	    \EE{\lnr \nabla f(X(t)\rnr_\bfr^2} &\leq 3\EE{\lnr \nabla f(X(t)) - \nabla f(\ones \mx(t)  )\rnr_\bfr^2} +  3\EE{\lnr \nabla f(\ones \mx(t))  - \ones \nabla f(\mx(t)) \rnr_\bfr^2}\nonumber\\
	    &\phantom{=} + 3\EE{ \lnr \ones \nabla f(\mx(t)) \rnr_\bfr^2}. 
	    \end{align}
	    Next, we bound each term in~\eqref{eq:norm-grad-f-X}. Recall from~\eqref{eq:def:grad-f-X} that the $i$th row of $\nabla f(X(t))$ is given by $\nabla f_i(X_i(t))$, where $X_i(t) = \bfx_i(t)$ is the $i$th row of matrix $X(t)$. Thus, from Assumption~\ref{asm:f}-(a) we have
	    \begin{align}
	    \EE{\lnr \nabla f(X(t)) - \nabla f(\ones \mx(t) )\rnr_\bfr^2} &= \EE{\sum_{i=1}^n r_i  \lnr \nabla f_i(\bfx_i(t)) - \nabla f_i( \mx(t)  )\rnr^2}  \cr
	    &\leq \li^2 \EE{\sum_{i=1}^n r_i \lnr \bfx_i(t)  -   \mx(t)  \rnr^2}  \cr
	    &= \li^2 \EE{\lnr X(t)  -  \ones  \mx(t) \rnr_\bfr^2}. 
	    \label{eq:norm-grad-f-X:1}
	\end{align}
    Similarly, using the convexity of function $\lnr \cdot \rnr^2$, for the second term in~\eqref{eq:norm-grad-f-X} we have 
    
    \begin{align}
        &\EE{\lnr \nabla f(\ones \mx(t))  - \ones \nabla f(\mx(t)) \rnr_\bfr^2}\cr
        &= \EE{\sum_{i=1}^n r_i \lnr \nabla f_i(\mx(t)) -\nabla  f(\mx(t))\rnr^2}\cr
        &= \sum_{i=1}^n r_i \EE{ 4 \lnr  \frac{1}{2}\Big(\nabla f_i(\mx(t)) -\nabla L(\mx(t)) \Big) - \frac{1}{2}\Big(   \nabla  f(\mx(t)) - \nabla L(\mx(t)) \Big)  \rnr^2}\cr
        &\leq \sum_{i=1}^n r_i \EE{2 \lnr \nabla f_i(\mx(t)) -\nabla L(\mx(t)) \rnr^2 + 2\lnr   \nabla  f(\mx(t)) - \nabla L(\mx(t))\rnr^2}\cr
        &= \sum_{i=1}^n 2r_i \EE{ \lnr \nabla f_i(\mx(t)) -\nabla L(\mx(t)) \rnr^2 } + 2\EE{ \lnr  \nabla  f(\mx(t)) - \nabla L(\mx(t)) \rnr^2}\cr
        &\stackrel{\rm{(a)}}{=} \!\sum_{i=1}^n 2r_i \mathbb{E}\!\lb\!\frac{1}{m_i^2}\! \lnr \sum_{j=1}^{m_i} \!\lb\nabla \ell(\mx(t),\xi_j^i)\! -\!\nabla L(\mx(t))\rb \rnr^2\rb \!\!+\!2\mathbb{E}\!\lb\! \frac{1}{N^2}\!\lnr  \sum_{j=1}^N \!\lb \nabla\ell(\mx(t), \xi_j) \!-\! \nabla L(\mx(t)) \rb \rnr^2 \rb\cr
        &\stackrel{\rm{(b)}}{=}  \sum_{i=1}^n 2r_i \frac{1}{m_i^2}  \sum_{j=1}^{m_i}  \EE{\lnr \nabla \ell(\mx(t),\xi_j^i) -\nabla L(\mx(t))\rnr^2}  + \frac{2}{N^2} \sum_{j=1}^N \EE{\lnr  \nabla\ell(\mx(t), \xi_j) - \nabla L(\mx(t)) \rnr^2 }\cr
        &\stackrel{\rm{(c)}}{\leq} \sum_{i=1}^n 2\frac{m_i}{N} \frac{1}{m_i^2} m_i \sgm^2  + \frac{2}{N^2}N\sgm^2 \cr
        &= \frac{2(n+1)}{N} \sgm^2,
        \label{eq:norm-grad-f-X:2}
    \end{align}
    where in {\rm (a)} we replaced the definitions of $f_i(\mx(t))$ and $f(\mx(t))$ from~\eqref{eq:ERM2} and~\eqref{eq:ERM}, respectively, the equality in {\rm (b)} holds since $\xi_j$s are independent samples from the underlying distribution, and  {\rm (c)} follows from Assumption~\ref{asm:f}-(b) and the fact that $r_i=m_i/N$ for $i\in [n]$. 
    Finally, for the third term in~\eqref{eq:norm-grad-f-X}, we have
    \begin{align}
       \EE{ \lnr \ones \nabla f(\mx(t)) \rnr_\bfr^2} =  \EE{\sum_{i=1}^n r_i \lnr \nabla f(\mx(t)) \rnr^2} = \EE{\lnr \nabla f(\mx(t)) \rnr^2}.
       \label{eq:norm-grad-f-X:3}
    \end{align}
     Plugging
     \eqref{eq:norm-grad-f-X:1}--\eqref{eq:norm-grad-f-X:3} in~\eqref{eq:norm-grad-f-X}, we get
    \begin{align}\label{eq:exp-norm-nabla-simple}
	    \EE{\lnr \nabla f(X(t)\rnr_\bfr^2}  &\leq 3\li^2 \EE{\lnr X(t)  -  \ones  \mx(t) \rnr_\bfr^2} +  3\frac{2(n+1)}{N} \sigma^2
	    +3 \EE{ \lnr  \nabla f(\mx(t)) \rnr^2}.
	\end{align}
    Next, replacing this bound in~\eqref{eq:total-dev}, we arrive at
	\begin{align}\label{eq:total-dev-bound} 
	\sum_{t=1}^{T} & \alpha(t)\beta(t)  \EE{\lnr X(t) -\ones \mx(t) \rnr_\bfr^2} \cr
	&\leq 
	    \frac{2\sm\kap}{\lambda} \sum_{t=1}^T \alpha(t) \beta^2(t) 
	     + 
	    \frac{8\kap}{\lambda^2}  \sum_{t=1}^T \lb \alpha^3(t) \beta(t)  \EE{ \Nr{\nabla f(X(t))}^2}\rb\cr
	&\leq 
	    \frac{2\sm\kap}{\lambda} \sum_{t=1}^T \alpha(t) \beta^2(t) 
	    + \frac{24\kap\li^2}{\lambda^2} \sum_{t=1}^{T} \alpha^3(t) \beta(t) \EE{\lnr X(t)  -  \ones  \mx(t) \rnr_\bfr^2} \cr
	    &\phantom{\leq} +  \frac{48\kap(n+1)\sgm^2}{N\lambda^2} \sum_{t=1}^{T} \alpha^3(t) \beta(t) 
	      +  \frac{24\kap}{\lambda^2}\sum_{t=1}^{T} \alpha^3(t) \beta(t)\EE{ \lnr  \nabla f(\mx(t)) \rnr^2}.
	\end{align}
	Now, let us define $\phi_{i,j}(T) := \sum_{t=1}^T \alpha^i(t) \beta^j(t)$.\\
	Then, ${\frac{2\sm\kap}{\lambda} \sum_{t=1}^T \alpha(t)    \beta^2(t) = \OneTwo \phi_{1,2}(T)}$ and ${\frac{48\kap(n+1)\sgm^2}{N\lambda^2} \sum_{t=1}^T \alpha^3(t) \beta(t) = \ThreeOne \phi_{3,1}(T)}$, where $\OneTwo := \frac{2\sm\kap}{\lambda}$ and $\ThreeOne := \frac{48\kap(n+1)\sgm^2}{N\lambda^2}$. 
	Furthermore, we 
	set $T_0:= \left\lceil \lp\frac{14 \azr \kap^{\frac{1}{2}}\li }{\lambda}\rp^{\frac{1}{\nu}}\right\rceil $  such that $\frac{24\kap\li^2}{\lambda^2}  \alpha^2(T_0) \leq \frac{24}{196}< \frac{1}{2} $.  Then, for $T > T_0$ we can rewrite~\eqref{eq:total-dev-bound} as
	\begin{align*}
	    \sum_{t=1}^{T_0} \alpha(t)\beta(t)& \EE{\Nr{X(t)-\ones \bar{\bfx}(t)}^2} + \sum_{t=T_0+1}^{T} \alpha(t)\beta(t) \EE{\Nr{X(t)-\ones\bar{\bfx}(t)}^2}\cr
     &\leq 
	     \OneTwo \phi_{1,2}(T) + \ThreeOne \phi_{3,1}(T) + \frac{24\kap}{\lambda^2}\sum_{t=1}^{T} \alpha^3(t) \beta(t)\EE{ \lnr  \nabla f(\mx(t)) \rnr^2}\cr 
	    &\phantom{=} 
	      + \frac{24\kap \li^2}{\lambda^2}  \sum_{t=1}^{T_0} \alpha^3(t) \beta(t) \EE{\lnr X(t)  -  \ones \mx(t) \rnr_\bfr^2} \cr 
	    &\phantom{=} 
	    +\frac{24 \kap \li^2}{\lambda^2} \sum_{t=T_0+1}^{T} \alpha^3(t) \beta(t) \EE{\lnr X(t)  -  \ones \mx(t) \rnr_\bfr^2}\cr
	 &\leq      
	      \OneTwo \phi_{1,2}(T) + \ThreeOne \phi_{3,1}(T) + \frac{24\kap}{\lambda^2}\sum_{t=1}^{T} \alpha^3(t) \beta(t)\EE{ \lnr  \nabla f(\mx(t)) \rnr^2}\cr 
	    &\phantom{=} 
	      + \frac{24 \kap \li^2}{\lambda^2}  \sum_{t=1}^{T_0} \alpha^3(t) \beta(t) \EE{\lnr X(t)  -  \ones \mx(t) \rnr_\bfr^2}\cr
	     &\phantom{=} + \frac{24 \kap \li^2}{\lambda^2} \alpha^2(T_0) \sum_{t=T_0+1}^{T} \alpha(t) \beta(t) \EE{\lnr X(t)  -   \ones\mx(t)  \rnr_\bfr^2}\cr
     &\leq      
	      \OneTwo \phi_{1,2}(T) + \ThreeOne \phi_{3,1}(T) + \frac{24 \kap}{\lambda^2}\sum_{t=1}^{T} \alpha^3(t) \beta(t)\EE{ \lnr  \nabla f(\mx(t)) \rnr^2}\cr 
	    &\phantom{=} 
	      + \SumTOne
	      + \frac{1}{2} \sum_{t=T_0+1}^{T} \alpha(t) \beta(t) \EE{\lnr X(t)  -  \ones \mx(t)  \rnr_\bfr^2},  
	\end{align*}
	where the second inequality holds since  $\alpha(t)$ is a non-increasing sequence, and
	\begin{align}
    	\SumTOne := \frac{24\kap \li^2}{\lambda^2} \sum_{t=1}^{T_0} \alpha^3(t) \beta(t) \EE{\lnr X(t)  -  \ones \mx(t) \rnr_\bfr^2 }
    	\label{eq:sum-T1}
	\end{align}
	does not grow with $T$. Therefore, we have
		\begin{align}\label{eq:total-dev-bound-T1}
		&\sum_{t=1}^{T} \alpha(t)\beta(t) \EE{\Nr{X(t)-\ones\bar{\bfx}(t)}^2}\cr 
		&\leq 
		2 \sum_{t=1}^{T_0} \alpha(t)\beta(t) \EE{\Nr{X(t)-\ones \bar{\bfx}(t)}^2} + \sum_{t=T_0+1}^{T} \alpha(t)\beta(t) \EE{\Nr{X(t)-\ones\bar{\bfx}(t)}^2} \cr
	    &\leq
	     2\OneTwo \phi_{1,2}(T) +2 \ThreeOne \phi_{3,1}(T) + 2 \SumTOne + \frac{48\kap}{\lambda^2}
	     \sum_{t=1}^{T} \alpha^3(t) \beta(t)\EE{ \lnr  \nabla f(\mx(t)) \rnr^2}.
	\end{align}
	
\subsection{Back to the Main Dynamics}	

Recall the dynamics in~\eqref{eq:main-ineq}, that is,
	\begin{align*}
     \EE{f\left(\mx(t+1)\right)}
         &\leq 
        \EE{ f\lp \mx(t)\rp} -\frac{1}{2}\at\et\EE{\lnr\nabla f(\mx(t))\rnr^2}\cr
        &\quad -\frac{1}{2} \at\et\lp 1-\at\et \li\rp\EE{\lnr \bfr^T \nabla f(X(t))\rnr^2} \cr
        &\quad + \at\et \li^2\EE{\Nr{X(t)  - \ones \mx(t)}^2} 
        + \beta^2(t)\frac{\li}{2}\sm.
    \end{align*}
Summing~\eqref{eq:main-ineq} for $t=1,2,\dots, T$, and using~\eqref{eq:total-dev-bound-T1} we get 
	\begin{align}\label{eq:sum-updates:1}
        &\EE{f\left(\mx(T+1)\right)}\cr
    &\leq 
        \EE{ f\lp \mx(1)\rp} -\frac{1}{2}\sum_{t=1}^{T}\at\et\EE{\lnr\nabla f(\mx(t))\rnr^2}\cr
        &\quad -\frac{1}{2}\sum_{t=1}^{T} \at\et\lp 1-\at\et \li\rp\EE{\lnr \bfr^T \nabla f(X(t))\rnr^2}  \cr 
        &\quad + \li^2 \sum_{t=1}^{T}\at\et \EE{\Nr{X(t)  -  \ones \mx(t)}^2} + \frac{\li}{2}\sm\sum_{t=1}^{T}\beta^2(t) \cr
    &\leq 
        \EE{ f\lp \mx(1)\rp} -\frac{1}{2}\sum_{t=1}^{T}\at\et\EE{\lnr\nabla f(\mx(t))\rnr^2} +  \frac{\li}{2}\sm \phi_{0,2} (T) \cr 
        &\quad -\frac{1}{2}\sum_{t=1}^{T} \at\et\lp 1-\at\et \li\rp\EE{\lnr \bfr^T \nabla f(X(t))\rnr^2}\cr
        &\quad + \li^2\lb  2\OneTwo \phi_{1,2}(T) +2 \ThreeOne \phi_{3,1}(T) + 2 \SumTOne + \frac{48\kap}{\lambda^2}
	     \sum_{t=1}^{T} \alpha^3(t) \beta(t)\EE{ \lnr  \nabla f(\mx(t)) \rnr^2} \rb, 
    \end{align}
    where $ \phi_{0,2} (T) = \sum_{t=1}^T \beta^2(t)$. 
Recall that for $T_0 \hspace{0pt}=\hspace{0pt} \left\lceil \hspace{0pt}\lp\frac{14 \azr \kap^{\frac{1}{2}}\li }{\lambda}\hspace{0pt}\rp^{\frac{1}{\nu}}\hspace{0pt}\right\rceil$. Moreover, for $\lambda$ and $\kap$ defined in Lemma~\ref{lemma:transition2} we have $\lambda\leq 1$ and $\kap>1$. Therefore, for $t>T_0$ we have \[\alpha(t)\beta(t)K \hspace{-1pt}\leq \alpha(T_0) \beta(T_0)K \leq \alpha(T_0)K \leq \frac{\lambda}{14 \kap^{\frac{1}{2}} }<1.\]
 This allows us to lower  bound the second summation in~\eqref{eq:sum-updates:1} for $T> T_0$ as
\begin{align}\label{eq:sum-updates:Neg}
   \frac{1}{2}&\sum_{t=1}^{T} \at\et\lp 1-\at\et \li\rp\EE{\lnr \bfr^T \nabla f(X(t))\rnr^2} \cr
   &= \frac{1}{2} \sum_{t=1}^{T_0} \at\et\lp 1-\at\et \li\rp\EE{\lnr \bfr^T \nabla f(X(t))\rnr^2}\cr
   &\quad + 
   \frac{1}{2}\sum_{t=T_0+1}^{T} \at\et\lp 1-\at\et \li\rp\EE{\lnr \bfr^T \nabla f(X(t))\rnr^2} \cr
   &\geq \frac{1}{2} \sum_{t=1}^{T_0} \at\et\lp 1-\at\et \li\rp\EE{\lnr \bfr^T \nabla f(X(t))\rnr^2} := \CNeg,
\end{align}
where $\CNeg$ does not grow with $T$, and the inequality holds since $1-\alpha(t)\beta(t)K \geq 0$ for $t>T_0$. 
Similarly, for $T_0 = \left\lceil \lp\frac{14 \azr \kap^{\frac{1}{2}}\li }{\lambda}\rp^{\frac{1}{\nu}}\right\rceil$, we have $\frac{48\kap\li^2}{\lambda^2}  \alpha^2(T_0) \leq \frac{48}{196}\azr \leq \frac{1}{4} $. Therefore, for $T > T_0$, the last summation in~\eqref{eq:sum-updates:1} can be upper bounded by
 	\begin{align}\label{eq:last-term-in-eq:sum-updates:1} 
           &\frac{48\kap \li^2}{\lambda^2}\sum_{t=1}^{T} \alpha^3(t) \beta(t)\EE{ \lnr  \nabla f(\mx(t)) \rnr^2}\cr
    &= 
         \frac{48\kap \li^2}{\lambda^2}  \sum_{t=1}^{T_0} \alpha^3(t) \beta(t)\EE{ \lnr  \nabla f(\mx(t)) \rnr^2} 
        + \frac{48\kap \li^2}{\lambda^2}  \sum_{t=T_0+1}^{T} \alpha^3(t) \beta(t)\EE{ \lnr  \nabla f(\mx(t)) \rnr^2}\cr
    &\leq 
        \SumTTwo
        + \frac{48\kap \li^2}{\lambda^2} \alpha^2(T_0)  \sum_{t=T_0+1}^{T} \alpha(t) \beta(t)\EE{ \lnr  \nabla f(\mx(t)) \rnr^2}\cr
    &\leq
        \SumTTwo
        + \frac{1}{4}  \sum_{t=T_0+1}^{T} \alpha(t) \beta(t)\EE{ \lnr  \nabla f(\mx(t)) \rnr^2},
    \end{align}
    where 
    \begin{align}
        \SumTTwo:= \frac{48\kap \li^2}{\lambda}  \sum_{t=1}^{T_0} \alpha^3(t) \beta(t)\EE{ \lnr  \nabla f(\mx(t)) \rnr^2},
        \label{eq:sum_T_2}
    \end{align}
    does not depend on $T$. Plugging~\eqref{eq:sum-updates:Neg} and \eqref{eq:last-term-in-eq:sum-updates:1} in \eqref{eq:sum-updates:1}, for $T>T_0$, we get 
    \begin{align*}
        \EE{ f\lp \mx(T+1)\rp}    
    &\leq 
        \EE{ f\lp \mx(1)\rp} + \frac{\li}{2} \phi_{0,2} (T)  + 2\li^2 \lp  \OneTwo \phi_{1,2}(T) + \ThreeOne \phi_{3,1}(T) + \SumTOne \rp \cr
        &\phantom{=}  
        -\frac{1}{2}\sum_{t=1}^{T} \at\et\lp 1-\at\et \li\rp\EE{\lnr \bfr^T \nabla f(X(t))\rnr^2}\cr
        &\phantom{=}      
        -\frac{1}{2}\sum_{t=1}^{T}\at\et\EE{\lnr\nabla f(\mx(t))\rnr^2} 
        + \frac{48\kap \li^2}{\lambda^2}\sum_{t=1}^{T} \alpha^3(t) \beta(t)\EE{ \lnr  \nabla f(\mx(t)) \rnr^2}\cr
    &\leq 
        \EE{ f\lp \mx(1)\rp} + \frac{\li}{2} \phi_{0,2} (T)  + 2\li^2 \lp  \OneTwo \phi_{1,2}(T) + \ThreeOne \phi_{3,1}(T) + \SumTOne\rp -\CNeg \cr
        &\phantom{=}      
        -\frac{1}{2}\sum_{t=1}^{T_0}\at\et\EE{\lnr\nabla f(\mx(t))\rnr^2}
        -\frac{1}{2}\sum_{t=T_0+1}^{T}\at\et\EE{\lnr\nabla f(\mx(t))\rnr^2} \cr
        &\phantom{=} + \SumTTwo
        + \frac{1}{4}  \sum_{t=T_0+1}^{T} \alpha(t) \beta(t)\EE{ \lnr  \nabla f(\mx(t)) \rnr^2}\cr
    &\leq 
        \EE{ f\lp \mx(1)\rp} + \frac{\li}{2} \phi_{0,2} (T)  + 2\li^2 \lp  \OneTwo \phi_{1,2}(T) + \ThreeOne \phi_{3,1}(T) + \SumTOne\rp-\CNeg + \SumTTwo\cr
        &\phantom{=}      
        -\frac{1}{4}\sum_{t=1}^{T}\at\et\EE{\lnr\nabla f(\mx(t))\rnr^2}.
    \end{align*}
    Hence, 
    \begin{align}
        &\sum_{t=1}^{T}\at\et\EE{\lnr\nabla f(\mx(t))\rnr^2} \cr
    &\leq 
        4 \EE{ f\lp \mx(1)\rp} - 4 \EE{f\left(\mx(T+1)\right)} 
        \!+\!  2\li \phi_{0,2} (T)  + 8\li^2 \lp  \OneTwo \phi_{1,2}(T) + \ThreeOne \phi_{3,1}(T) + \SumTOne\rp \!-\!4\CNeg \!+\! 4\SumTTwo
        \cr    
    &\leq 
        4 \EE{ f\lp \mx(1)\rp} - 4 \EE{f(\bfx^\star)} + 2\li \phi_{0,2} (T)  + 8\li^2 \lp\OneTwo \phi_{1,2}(T) + \ThreeOne \phi_{3,1}(T) + \SumTOne \rp -4\CNeg + 4\SumTTwo
        \cr
    &= \const + 2\li \phi_{0,2} (T) + 8\li^2 \OneTwo \phi_{1,2}(T) + 8\li^2 \ThreeOne \phi_{3,1}(T)\cr
    &\leq \const + (2\li + 8\li^2  \OneTwo \alpha_0) \phi_{0,2} (T) + 8\li^2 \ThreeOne \phi_{3,1}(T),
        \label{eq:main-bound-T}
    \end{align}
    where 
    \begin{align*}
        \const &:= 8  \li^2  \SumTOne -4\CNeg + 4 \SumTTwo + 4 \EE{ f\lp \mx(1)\rp} - 4 \EE{f(\bfx^\star)} 
         \cr
         &= \frac{192\kap \li^2}{\lambda}  \sum_{t=1}^{T_0} \alpha^3(t) \beta(t) \lb \li^2  \EE{\lnr X(t)  -  \ones \mx(t) \rnr_\bfr^2} + \EE{ \lnr  \nabla f(\mx(t)) \rnr^2} \rb \cr
         & \phantom{=} -2\sum_{t=1}^{T_0} \at\et\lp 1-\at\et \li\rp\EE{\lnr \bfr^T \nabla f(X(t))\rnr^2} + 4\lp \EE{ f\lp \mx(1)\rp} -  \EE{f(\bfx^\star)}\rp
    \end{align*}
    is constant (does not depend on $T$), and the last inequality in~\eqref{eq:main-bound-T} follows from the fact that 
    \begin{align}\label{eq:phi012}
     \phi_{1,2}(T) = \sum_{t=1}^T \alpha(t) \beta^2(t)   = \alpha_0 \sum_{t=1}^T \frac{1}{(t+\tau)^\nu} \beta^2(t) \leq \alpha_0 \sum_{t=1}^T  \beta^2(t) = \alpha_0 \phi_{0,2}(T).
    \end{align}
\subsection{Bound on the Moments of \texorpdfstring{$\EE{\lnr\nabla f(\mx(t))\rnr^2}$}{} and \texorpdfstring{$\EE{\Nr{X(t)-\ones \bar{\bfx}(t) }^2}$}{}}
\begin{figure}
    \centering
    \includegraphics[width=0.5\textwidth]{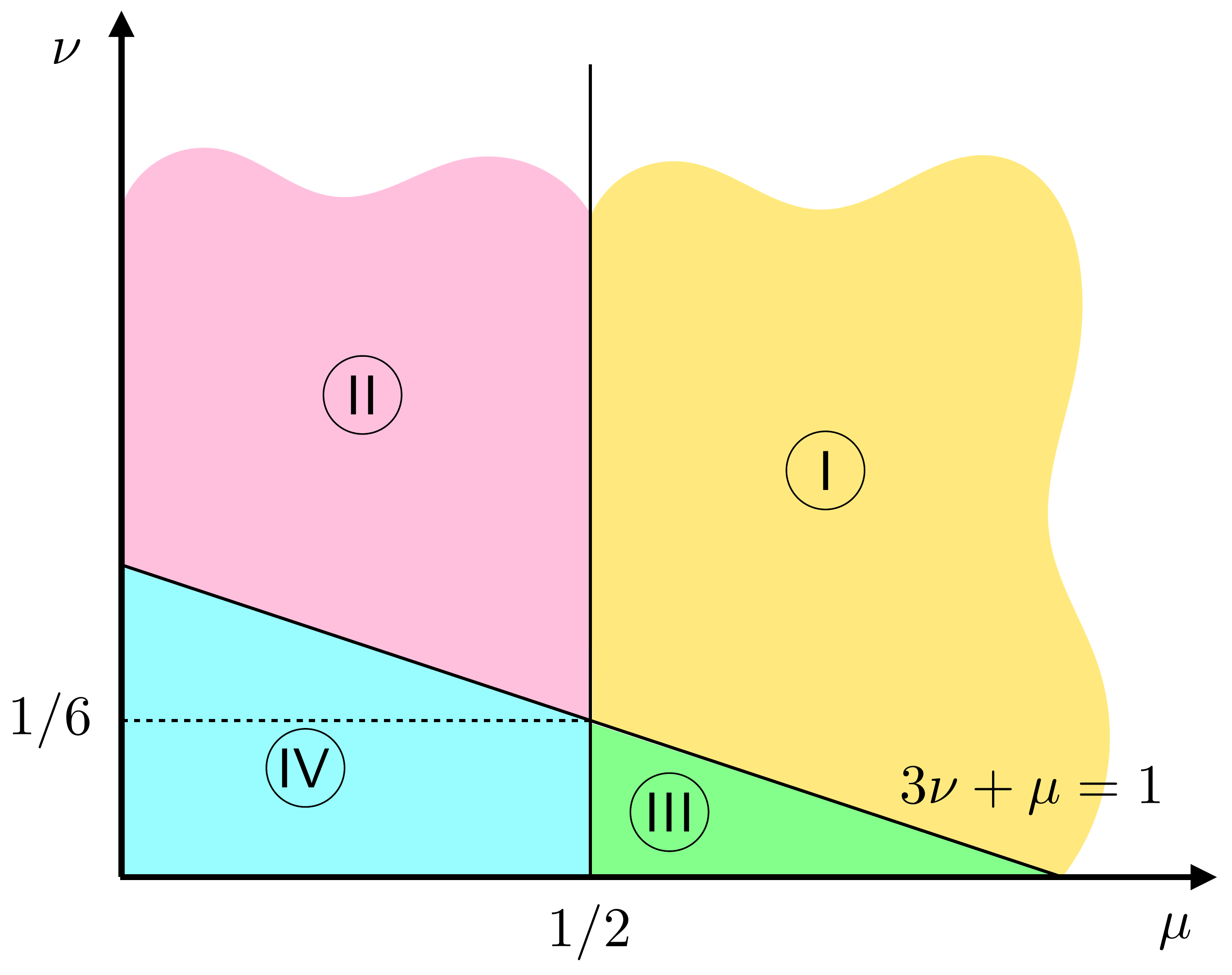}
    \caption{Regions of $(\mu,\nu)$.}
    \label{fig:regions}
\end{figure}

The inequality in~\eqref{eq:main-bound-T} provides us with an upper bound on the temporal average of sequence $\big\{ \at\et\mathbb{E}[\lnr\nabla f(\mx(t))\rnr^2]\big\}$. However, our goal is to derive a bound on the temporal average of $\big\{ \mathbb{E}[\lnr\nabla f(\mx(t))\rnr^2]\big\}$. To this end,  for any given $\theta\in (0,1)$, we  define the measure of convergence as
\begin{align*}
    &M_\theta(\nu,\mu) := \lb \frac{1}{T} 
    \sum_{t=1}^T \lp\EE{\lnr\nabla f(\mx(t))\rnr^2}\rp^\theta \rb^{\frac{1}{\theta}}.
\end{align*}
Note that by H\"older's inequality \cite[Theorem\ 6.2]{folland1999real} for any $p,q>1$  with $\frac{1}{p}+\frac{1}{q}=1$, and non-negative sequences $\{a_t\}_{t=1}^T$ and $\{b_t\}_{t=1}^T$, we have 
\begin{align*}
    \sum_{t=1}^T a_t b_t \leq \lp \sum_{t=1}^T a_t^p \rp^{\frac{1}{p}} \lp \sum_{t=1}^T b_t^q \rp^{\frac{1}{q}},
\end{align*}
or equivalently,
\begin{align}\label{eq:holder}
    \lp \sum_{t=1}^T a_t b_t \rp^q \leq \lp \sum_{t=1}^T a_t^p \rp^{\frac{q}{p}} \lp \sum_{t=1}^T b_t^q \rp.
\end{align}
Let $(p,q)=\lp \frac{1}{1-\theta},\frac{1}{\theta}\rp$ so that $\frac{1}{p}+ \frac{1}{q}=1$. Furthermore, let  
\begin{align}\label{eq:at-bt}
\begin{split}
    a_t &= \left(\frac{1}{\alpha(t) \beta(t)}\right)^{\theta} = \frac{1}{(\azr\bzr)^{\theta}} (t+\tu)^{(\mu+\nu)\theta}\\
    b_t &=  \lp \at\et\EE{\lnr\nabla f(\mx(t))\rnr^2}\rp^\theta.
    \end{split}
\end{align}
Then, applying H\"older's inequality \eqref{eq:holder}, we arrive at 
\begin{align}\label{eq:holder-RHS}
M_\theta(\nu,\mu) &= \lb \frac{1}{T} 
    \sum_{t=1}^T \lp\EE{\lnr\nabla f(\mx(t))\rnr^2}\rp^\theta \rb^{\frac{1}{\theta}} = \lp \frac{1}{T}  \sum_{t=1}^T a_t b_t\rp^q
    \leq \frac{1}{T^q}  \lp \sum_{t=1}^T a_t^p \rp^{\frac{q}{p}}  \lp \sum_{t=1}^T b_t^q \rp. 
\end{align}
It remains to upper bound the terms in the right hand side of~\eqref{eq:holder-RHS}.
First, using Lemma~\ref{lm:sum-exp} we get 
\begin{align*}
    \sum_{t=1}^T a_t^p =  \frac{1}{(\alpha_0 \beta_0)^{\frac{\theta}{1-\theta}}} \sum_{t=1}^T (t+\tau)^{\frac{(\nu + \mu)\theta}{1-\theta}} \leq \frac{2^{1+\frac{(\nu + \mu)\theta}{1-\theta} }} {(\alpha_0 \beta_0)^{\frac{\theta}{1-\theta}}  \lp 1+ \frac{(\nu + \mu)\theta}{1-\theta} \rp } (T+\tau)^{1+\frac{(\nu + \mu)\theta}{1-\theta}}.
\end{align*}
Therefore, 
\begin{align}\label{eq:holder-B}
    \lp \sum_{t=1}^T a_t^p \rp^{\frac{q}{p}}  \leq \frac{2^{\frac{1-\theta}{\theta}+(\nu + \mu)}} {\alpha_0  \beta_0   \lp 1+ \frac{(\nu + \mu)\theta}{1-\theta} \rp^{\frac{1-\theta}{\theta}} } (T+\tau)^{\frac{1-\theta}{\theta}+\nu + \mu}:=\cntB (T+\tau)^{\frac{1-\theta}{\theta}+\nu + \mu}.
\end{align}
Next, continuing from~\eqref{eq:main-bound-T}, for $T\geq T_0$ we can write 
\begin{align}\label{eq:holder-A}
     \sum_{t=1}^T b_t^q&=  \sum_{t=1}^T \at\et\EE{\lnr\nabla f(\mx(t))\rnr^2} \cr
    & \leq  \const + (2\li + 8\li^2  \OneTwo \alpha_0) \phi_{0,2} (T) + 8\li^2 \ThreeOne \phi_{3,1}(T) \cr
    &\leq  \max\lc 3\const , (6\li + 24\li^2  \OneTwo \alpha_0) \phi_{0,2} (T) , 24\li^2 \ThreeOne  \phi_{3,1}(T) \rc,
\end{align}
where the last inequality follows from the fact that $a+b+c\leq \max\{3a,3b,3c\}$. 
Next, we can use  Lemma~\ref{lm:sum-exp} to bound $\phi_{0,2}(T)$, as
\begin{align}\label{eq:mu-reg}
    \phi_{0,2}(T) = \sum_{t=1}^T \beta^2(t) = \frac{1}{\beta_0^2} \sum_{t=1}^T (t+\tau)^{-2\mu} \leq 
    \begin{cases}
    \frac{\tau^{1-2\mu}}{\beta_0^2 |1-2\mu|} & \textrm{if $\mu>1/2$},\\
    \frac{1}{\beta_0^2}\ln\lp \frac{T}{\tau} +1 \rp& \textrm{if $\mu=1/2$},\\
    \frac{2}{\beta_0^2(1-2\mu)} (T+\tau)^{1-2\mu} & \textrm{if $0\leq \mu<1/2 $},
    \end{cases}
\end{align}
where we used the fact that $2^{1-2\mu}\leq 2$ for $0\leq \mu <1$. Similarly, applying Lemma~\ref{lm:sum-exp} on $\phi_{3,1}(T)$, we get
\begin{align}\label{eq:nu-mu-reg}
    \phi_{3,1}(T) &= \sum_{t=1}^T \alpha^3(t) \beta(t) = \frac{1}{\alpha_0^3 \beta_0} \sum_{t=1}^T (t+\tau)^{-3\nu-\mu} \cr
    &\leq 
    \begin{cases}
    \frac{\tau^{1-3\nu-\mu}}{\alpha_0^3 \beta_0 |1-3\nu-\mu|} & \textrm{if $3\nu+\mu>1$},\\
    \frac{1}{\alpha_0^3 \beta_0}\ln\lp \frac{T}{\tau} +1 \rp& \textrm{if $3\nu+\mu=1$},\\
    \frac{2}{\alpha_0^3 \beta_0(1-3\nu-\mu)} (T+\tau)^{1-3\nu-\mu} & \textrm{if $0\leq 3\nu+\mu<1 $},
    \end{cases}
\end{align}
in which we have upper bounded $2^{1-3\nu-\mu}$ by  by $2$, for $0\leq 3\nu+\mu<1 $. 
Note that the upper bound in~\eqref{eq:holder-A} is the maximum of a constant, and two $\phi_{\cdot,\cdot}(T)$ functions. Moreover, depending on the values of $\mu$ and $\nu$, each $\phi_{\cdot,\cdot}(T)$ can be upper bounded by a constant, a logarithmic function, or a polynomial (in $T$). Figure~\ref{fig:regions} illustrates the four regions of $(\mu,\nu)$. In the following we first, analyze the interior of the four regions shown in Figure~\ref{fig:regions}, and then study the boundary cases. 

\noindent \underline{\textbf{Region $\mathsf{I}$: $\mu>1/2$ and $3\nu+\mu >1$}}\\
Recall from~\eqref{eq:mu-reg} and ~\eqref{eq:nu-mu-reg} that both $\phi_{0,2}(T)$ and $\phi_{3,1}(T)$ are upper bounded by constant functions, which grows faster than a constant. Hence, from~\eqref{eq:holder-A} in this regime we have 
\begin{align}\label{eq:holder-A-I}
     \sum_{t=1}^T b_t^q 
    & \leq  \max \lc 3\const , (6\li + 24\li^2  \OneTwo \alpha_0) \phi_{0,2} (T) , 24\li^2 \ThreeOne \phi_{3,1}(T)\rc  \cr
    &\leq   \max \lc 3\const , (6\li + 24\li^2  \OneTwo \alpha_0) \frac{\tau^{1-2\mu}}{\beta_0^2 |1-2\mu|}  , 24\li^2 \ThreeOne  \frac{\tau^{1-3\nu -\mu}}{\alpha_0^3\beta_0 |1-3\nu -\mu|} \rc \cr
    &:=\COne.
\end{align}
Note that $\COne$ does not depend on $T$. Plugging \eqref{eq:holder-B} and~\eqref{eq:holder-A-I} in~\eqref{eq:holder-RHS}, we arrive at 
\begin{align}
M_\theta(\nu,\mu) 
    \leq \frac{1}{T^{1/\theta}}  \cntB \cdot \COne\cdot  (T+\tau)^{\frac{1-\theta}{\theta}+\nu+\mu} = \BO\lp T^{-(1-\nu-\mu)}\rp.
\end{align}

\noindent \underline{\textbf{Region $\mathsf{II}$: $\mu<1/2$ and $3\nu+\mu >1$}}\\ 
From~\eqref{eq:mu-reg} and~\eqref{eq:nu-mu-reg}, we can see that $\phi_{0,2}(T)$ and $\phi_{3,1}(T)$ the upper bounded by a polynomial and a constant function of $T$, respectively. Noting that for sufficiently large  $T$, a polynomial function beats any constant, we can write
\begin{align}\label{eq:holder-A-II}
    \sum_{t=1}^T b_t^q
    & \leq   \max \lc 3\const , (6\li + 24\li^2  \OneTwo \alpha_0) \phi_{0,2} (T) , 24\li^2 \ThreeOne \phi_{3,1}(T)\rc  \cr
    &\leq  \frac{12\li + 48\li^2  \OneTwo \alpha_0}{\beta_0^2 (1-2\mu)} (T+\tau)^{1-2\mu}\cr
    &:=\CTwo \cdot (T+\tau)^{1-2\mu}.
\end{align}
This together with~\eqref{eq:holder-RHS} and~\eqref{eq:holder-B} implies
\begin{align}
M_\theta(\nu,\mu) 
    \leq \frac{1}{T^{1/\theta}}  \cntB (T+\tau)^{\frac{1-\theta}{\theta}+\nu+\mu}
     \cdot\CTwo  \cdot (T+\tau)^{1-2\mu} = \BO\lp T^{-(\mu-\nu)}\rp.
\end{align}

\noindent \underline{\textbf{Region $\mathsf{III}$: $\mu>1/2$ and $3\nu+\mu < 1$}}\\
Recall from~\eqref{eq:mu-reg} and~\eqref{eq:nu-mu-reg} that  $\phi_{0,2}(T)$ and $\phi_{3,1}(T)$ are upper bounded by a constant and a polynomial function of $T$, respectively. Therefore, for sufficiently large $T$, we get 
\begin{align}\label{eq:holder-A-III}
     \sum_{t=1}^T b_t^q 
    & \leq  \max \lc 3\const , (6\li + 24\li^2  \OneTwo \alpha_0) \phi_{0,2} (T) , 24\li^2 \ThreeOne \phi_{3,1}(T)\rc \cr
    &\leq \lp  \frac{48\li^2 \ThreeOne }{\alpha_0^3 \beta_0 (1-3\nu-\mu)} \rp (T+\tau)^{1-3\nu-\mu}\cr
    &:=\CThree \cdot (T+\tau)^{1-3\nu-\mu}.
\end{align}
Therefore, plugging~\eqref{eq:holder-A-III} and~\eqref{eq:holder-B} to~\eqref{eq:holder-RHS}, we get 
\begin{align}
M_\theta(\nu,\mu) 
    \leq \frac{1}{T^{1/\theta}}  \cntB (T+\tau)^{\frac{1-\theta}{\theta}+\nu+\mu}
   \cdot  \CThree \cdot  (T+\tau)^{1-3\nu-\mu} = \BO\lp T^{-2\nu}\rp.
\end{align}

\noindent \underline{\textbf{Region $\mathsf{IV}$: $\mu<1/2$ and $3\nu+\mu < 1$}}\\
In this region, we can use~\eqref{eq:mu-reg} and~\eqref{eq:nu-mu-reg} to upper bound both  $\phi_{0,2}(T)$ and $\phi_{3,1}(T)$ by polynomial functions. Thus, for sufficiently large $T$, we get 
\begin{align}\label{eq:holder-A-IV}
    \sum_{t=1}^T b_t^q 
    & \leq  \max \lc 3\const , (6\li + 24\li^2  \OneTwo \alpha_0) \phi_{0,2} (T) , 24\li^2 \ThreeOne \phi_{3,1}(T)\rc \cr
    &\leq  \max \lc  \frac{12\li + 48\li^2  \OneTwo \alpha_0}{\beta_0^2 (1-2\mu)} ,  \frac{48\li^2 \ThreeOne}{\alpha_0^3 \beta_0 (1-3\nu-\mu)} \rc 
    (T+\tau)^{ \max\{1-2\mu, 1-3\nu-\mu\} }
    \cr
    &:=\CFour \cdot (T+\tau)^{\max\{1-2\mu, 1-3\nu-\mu\} }.
\end{align}
Then, we can plug~\eqref{eq:holder-A-IV} and~\eqref{eq:holder-B} to~\eqref{eq:holder-RHS}, to conclude
\begin{align}
M_\theta(\nu,\mu) 
    &\leq \frac{1}{T^{1/\theta}}  \cntB (T+\tau)^{\frac{1-\theta}{\theta}+\nu+\mu}\cdot \CFour \cdot (T+\tau)^{\max\{1-2\mu, 1-3\nu-\mu\}  } \cr
    &= \BO\lp T^{-\min\{\mu-\nu, 2\nu\}  }\rp.
\end{align}
Summarizing the result of the four cases above, we get 
\begin{align}
   M_\theta(\nu,\mu) &= \begin{cases}
   \BO\lp T^{-(1-\nu-\mu)}\rp & \textrm{if $\mu>1/2$ and $3\nu+\mu>1$}\\
   \BO\lp T^{-(\mu-\nu)}\rp & \textrm{if $\mu<1/2$ and $3\nu+\mu>1$}\\
   \BO\lp T^{-2\nu}\rp & \textrm{if $\mu>1/2$ and $3\nu+\mu<1$}\\
   \BO\lp T^{-\min\{\mu-\nu, 2\nu\} }\rp & \textrm{if $\mu<1/2$ and $3\nu+\mu<1$}
   \end{cases}\\
   &= \BO\lp T^{-\min\{1-\nu-\mu, \mu-\nu, 2\nu\} }\rp,
\end{align}
which concludes the first claim of Theorem~\ref{thm:non_cvx}. \\
Recall that our goal is to optimize over $(\nu,\mu)$ to achieve the best convergence rate. This is equivalent to maximizing the exponent of $1/T$ in each of the four regions, i.e., 
\begin{eqnarray*}
    &\textrm{Region }\mathsf{I}:\quad &\sup \{1-\nu-\mu: \mu>1/2, 3\nu+\mu>1\},\\
    &\textrm{Region }\mathsf{II}:\quad &\sup \{\mu-\nu: \mu<1/2, 3\nu+\mu>1\},\\
    &\textrm{Region }\mathsf{III}:\quad &\sup \{2\nu: \mu>1/2, 3\nu+\mu<1\},\\
    &\textrm{Region }\mathsf{IV}:\quad &\sup \{\min\{\mu-\nu, 2\nu\}: \mu<1/2, 3\nu+\mu<1\}.
\end{eqnarray*}
Interestingly, it turns out that the respective supremum value in all four regions is $1/3$, which corresponds to the boundary point $(\nu^\s,\mu^\s)=(1/6,1/2)$. However,  this point does not belong to any of the corresponding open sets, which motivates the convergence analysis of $M_{\theta}$ for $(\nu^\s,\mu^\s)=(1/6,1/2)$. 

\noindent \underline{\textbf{Boundary Case:   $\mu=1/2$ and $3\nu+\mu = 1$}}\\
First note that the two lines of interest intersect at $(\nu^\s,\mu^\s)=(1/6,1/2)$, as shown in Figure~\ref{fig:regions}.   Applying Lemma~\ref{lm:sum-exp} on $\phi_{0,2}(T)$ and $\phi_{3,1}(T)$ for $(\nu^\s,\mu^\s)=(1/6,1/2)$ we get
\begin{align}
    \phi_{0,2}(T)  &\leq \frac{1}{\beta_0^2}\ln\lp \frac{T}{\tau}+1\rp,\cr
    \phi_{3,1}(T) &\leq \frac{1}{\alpha_0^3\beta_0} \ln \lp\frac{T}{\tau}+1\rp.
\end{align}
Therefore, \eqref{eq:holder-A} reduces to 
\begin{align}\label{eq:holder-A-Bndr}
   \sum_{t=1}^T b_t^q 
    &\leq 
     \max\lc 3\const , (6\li + 24\li^2  \OneTwo \alpha_0) \phi_{0,2} (T) , 24\li^2 \ThreeOne  \phi_{3,1}(T) \rc \cr
    &\leq  \max \lc  \frac{6\li + 24\li^2  \OneTwo \alpha_0}{\beta_0^2}  , \frac{24\li^2 \ThreeOne}{\alpha_0^3 \beta_0}  \rc  \ln \lp\frac{T}{\tau}+1\rp  \cr
    &:= \CBndr \cdot \ln \lp\frac{T}{\tau}+1\rp .
\end{align}
Plugging~\eqref{eq:holder-A-Bndr} and~\eqref{eq:holder-B} to~\eqref{eq:holder-RHS}, we arrive at 
\begin{align}
M_\theta\lp \frac{1}{6},\frac{1}{2}\rp
    &\leq \frac{1}{T^{1/\theta}}  \cntB (T+\tau)^{\frac{1-\theta}{\theta}+\frac{2}{3}} \cdot\CBndr \cdot \ln \lp\frac{T}{\tau}+1\rp \cr
    &= \BO\lp T^{-1/3} \ln T\rp,
\end{align}
which is the second claim of the theorem. 
\section{Proof of Proposition~\ref{cor:M1}}
\label{sec:prf:cor}
First note that we cannot directly conclude the proposition from Theorem~\ref{thm:non_cvx}, since the theorem only holds for $\theta\in(0,1)$, and not $\theta=1$. In order to show the claim, for a vector $\bfy\in \mathbb{R}^T$ and some $\theta\in(0,1)$ we define\footnote{Note that for $\|\bfy\|_\theta$ is not a norm, since it is not a subadditive function for $\theta<1$.} 
$\|\bfy\|_\theta:=\lp|y_1|^\theta+|y_2|^\theta+\ldots+ |y_T|^\theta\rp^{1/\theta}
$. 
Then we have $\|\bfy\|_1 \leq \|\bfy\|_\theta$ since
\begin{align*}
    \lp \frac{\|\bfy\|_\theta}{\|\bfy\|_1}\rp^\theta = \frac{|y_1|^\theta+\ldots+ |y_T|^\theta}{(|y_1|+\ldots + |y_T|)^\theta} = \sum_{t=1}^T \lp \frac{|y_t|}{|y_1|+\ldots + |y_T|} \rp^\theta \geq \sum_{t=1}^T  \frac{|y_t|}{|y_1|+\ldots +|y_T|} =1,
\end{align*}
where the inequality holds since we have $0\leq |y_t|/\sum_{i=1}^T |y_i| \leq 1$, and  $\theta<1$. \\
Now, for the vector $\bfy$ with $y_t=\EE{\lnr\nabla f(\mx(t))\rnr^2}$ we have 
\begin{align*}
    M_1 &\!=\! \frac{1}{T} \sum_{t=1}^T \EE{\lnr\nabla f(\mx(t))\rnr^2} \!=\! \frac{1}{T} \|\bfy\|_1 \!\leq\! \frac{1}{T} \|\bfy\|_\theta = \frac{1}{T}   \lp \sum_{t=1}^T \lp \EE{\lnr\nabla f(\mx(t))\rnr^2}\rp^\theta \rp^{1/\theta} \\
    &= \frac{1}{T^{1-1/\theta}} \lp \frac{1}{T}\sum_{t=1}^T \lp \EE{\lnr\nabla f(\mx(t))\rnr^2}\rp^\theta \rp^{1/\theta}
    = \frac{1}{T^{1-1/\theta}} M_{\theta}. 
\end{align*}
Then, from Theorem~\ref{thm:non_cvx} for $(\nu^\s,\mu^\s)=(1/6,1/2)$ and
$\theta = \frac{2}{2+\epsilon}$, we get 
\begin{align}
    M_1 \leq \frac{1}{T^{1-1/\theta}} M_{\theta} \leq \frac{1}{T^{1-1/\theta}} \BO \lp T^{-1/3} \ln T\rp = T^{\epsilon/2} \BO \lp T^{-1/3} \ln T\rp = \BO \lp T^{-1/3+\epsilon} \rp,
\end{align}
where the last equality holds since 
$\ln T = \BO(T^{\epsilon/2})$  for any $\epsilon>0$. Similarly, for the vector $\bfz\in\mathbb{R}^T$  with $z_t := \EE{\Nr{X(t) - \ones \mx(t) }^2}$ and any $\theta\in(0,1)$ we have
\begin{align}\label{eq:prf:cor}
    \frac{1}{T}\sum_{t=1}^T \EE{\Nr{X(t) - \ones \mx(t) }^2} = \frac{1}{T} \|\bfz\|_1  \leq \frac{1}{T} \|\bfz\|_\theta  &= \frac{1}{T}  \lb
    \sum_{t=1}^T \lp\EE{\Nr{X(t)\!-\!\ones \bar{\bfx}(t) }^2}\rp^\theta \rb^{\frac{1}{\theta}} 
    \nonumber\\
    &
    = \!\frac{1}{T^{1-1/\theta}}  \lb \frac{1}{T}
    \sum_{t=1}^T \lp\EE{\Nr{X(t)\!-\!\ones \bar{\bfx}(t) }^2}\rp^\theta \rb^{\frac{1}{\theta}}. 
\end{align}
We need to bound the RHS of~\eqref{eq:prf:cor}. Let $a_t\!:=\! (\alpha(t) \beta(t))^{-\theta}\! = (t+\tau)^{2\theta/3}/(\alpha_0 \beta_0)^\theta$  and 
$c_t := \lp\alpha(t)\beta(t) \EE{\Nr{X(t)-\ones \bar{\bfx}(t) }^2}\rp^{\theta}
$. Then,  using the H\"older's inequality in~\eqref{eq:holder} for $(p,q) =\lp \frac{1}{1-\theta}, \frac{1}{\theta}\rp$ we can write
\begin{align}\label{eq:holder-var}
 \lb \frac{1}{T} 
    \sum_{t=1}^T \lp\EE{\Nr{X(t)-\ones \bar{\bfx}(t) }^2}\rp^\theta \rb^{\frac{1}{\theta}} = \lp \frac{1}{T}  \sum_{t=1}^T a_t c_t\rp^q
    \leq \frac{1}{T^q}  \lp \sum_{t=1}^T a_t^p \rp^{\frac{q}{p}}  \lp \sum_{t=1}^T c_t^q \rp. 
\end{align}
Note that $\lp \sum_{t=1}^T a_t^p \rp^{\frac{q}{p}}$ is bounded in~\eqref{eq:holder-B}.
Moreover, for $\sum_{t=1}^T c_t^q$  we can write
\begin{align}
	&\sum_{t=1}^T c_t^q = \sum_{t=1}^{T} \alpha(t)\beta(t) \EE{\Nr{X(t)-\ones \bar{\bfx}(t) }^2}\nonumber\\ 
	 &\stackrel{\rm{(a)}}{\leq}
	 2\OneTwo \phi_{1,2}(T) +2 \ThreeOne \phi_{3,1}(T) + 2 \SumTOne + \frac{48\kap}{\lambda^2}
	 \sum_{t=1}^{T} \alpha^3(t) \beta(t)\EE{ \lnr  \nabla f(\mx(t))
	 \rnr^2}\nonumber\\
	 &\stackrel{\rm{(b)}}{\leq}
	  2\OneTwo \phi_{1,2}(T) +2 \ThreeOne \phi_{3,1}(T) + 2 \SumTOne +  \frac{48\kap \alpha_0^2}{\lambda^2}
	 \sum_{t=1}^{T} \alpha(t) \beta(t)\EE{ \lnr  \nabla f(\mx(t))
	 \rnr^2}\nonumber\\
	 &\stackrel{\rm{(c)}}{\leq}\!\!  2 \alpha_0 \OneTwo \phi_{0,2}(T) \!+\!2 \ThreeOne \phi_{3,1}(T) \!+\! 2 \SumTOne \!+\! \frac{48\kap \azr^2}{\lambda^2}\! \lp\! \const\! +\! \li(2\!+\! 8\li \! \alpha_0 \OneTwo ) \phi_{0,2} (T) \!+\! 8\li^2 \ThreeOne \phi_{3,1}(T) \rp\nonumber\\
	 & \leq \! 2\SumTOne\! +\! \frac{48\kap \azr^2\const}{\lambda^2} \!+\! \lp\! 2\azr\OneTwo \!+\! \frac{96\kap\alpha_0^2 K (1 + 4\li \alpha_0 \OneTwo)}{\lambda^2} \!\rp \phi_{0,2}(T) \!+\! \lp 2\ThreeOne \!+\!  \frac{384\kap \azr^2\li^2\ThreeOne }{\lambda^2}\!\rp\! \phi_{3,1}(T)\nonumber\\
	 &\stackrel{\rm{(d)}}{\leq} 2\SumTOne  + \frac{48\kap \azr^2\const}{\lambda^2} + \lp \frac{2\azr\OneTwo}{\beta_0^2} + \frac{96\kap\alpha_0^2 K (1 + 4\li \alpha_0 \OneTwo)}{\lambda^2\beta_0^2} \rp \ln\lp\frac{T}{\tau}+ 1\rp 
	 \nonumber\\&\hspace{79pt} 
	 + \lp \frac{2\ThreeOne}{\alpha_0^3 \beta_0} +  \frac{384\kap \li^2 \ThreeOne }{\lambda^2 \alpha_0 \beta_0}\rp \ln\lp\frac{T}{\tau}+ 1\rp \stackrel{\rm{(e)}}{\leq} \Cvar \cdot \ln\lp\frac{T}{\tau}+ 1\rp, \label{eq:bnd-var-1} 
\end{align}
where $\rm{(a)}$ follows from~\eqref{eq:total-dev-bound-T1},  $\rm{(b)}$ holds since $\alpha^2(t)=\frac{\azr^2}{(t+\tau)^{1/3}} \leq \alpha_0^2$, the inequality in $\rm{(c)}$ follows from~\eqref{eq:main-bound-T} and~\eqref{eq:phi012}, we have used a bounds in~\eqref{eq:mu-reg} and~\eqref{eq:nu-mu-reg} for $(\nu^\s,\mu^\s)=(1/6,1/2)$ in $\rm{(d)}$, and $\rm{(e)}$ holds since the constant term $2\SumTOne  + 48\kap \azr^2\const/\lambda^2$ is upper bounded by $\ln(T/\tau +1)$ for large enough $T$. 
Plugging~\eqref{eq:holder-B} for $\nu^\star+\mu^\star=2/3$ and~\eqref{eq:bnd-var-1} into~\eqref{eq:holder-var} and~\eqref{eq:prf:cor}, and setting $\theta = \frac{2}{2+\epsilon}$ we arrive at
\begin{align}\label{eq:var-theta}
 \frac{1}{T}\sum_{t=1}^{T} \mathbb{E}\Big[&\|\bfx_i(t)-\bar{\bfx}(t)\|^2\Big]
 \leq \frac{1}{\rmin} \lb 
     \frac{1}{T}\sum_{t=1}^{T}\EE{\|X(t)-\ones \bar{\bfx}(t)\|_{\bfr}^2}\rb  \nonumber\\
    &\leq\frac{1}{\rmin} \frac{1}{T^{1-1/\theta}} \frac{1}{T^{1/\theta}}  \cntB (T+\tau)^{\frac{1-\theta}{\theta}+\frac{2}{3}} \cdot \Cvar\cdot  \ln\lp\frac{T}{\tau}+ 1\rp 
    = \BO \lp T^{-1/3 +\epsilon}\rp, 
\end{align}
where the last equality holds since  $\ln T = \BO(T^{\epsilon/2})$  for any $\epsilon>0$. Finally, combining~\eqref{eq:thm:M1},~\eqref{eq:thm:var}, and using Assumption~\ref{asm:f} and   Lemma~\ref{lm:nor_squr} (for $\omega=1$)  we have
\begin{align*}
    \frac{1}{T}&\sum_{t=1}^{T} \mathbb{E}\Big[\|\nabla f( \bfx_i(t))\|^2\Big] \leq 
    \frac{2}{T}\sum_{t=1}^{T} \left\{\mathbb{E}\Big[\|\nabla f( \bfx_i(t))-\nabla f(\mx(t)) \|^2\Big] + \mathbb{E}\Big[\|\nabla f(\mx(t))\|^2\Big]\right\} \nonumber\\
    &\leq 
    \frac{2}{T}\ \sum_{t=1}^{T}\left\{ \li^2\mathbb{E}\Big[\|\bfx_i(t)-\bar{\bfx}(t)\|^2\Big] +\mathbb{E}\Big[\|\nabla f(\mx(t))\|^2\Big]\right\}\leq \BO(T^{-1/3 + \epsilon}).
\end{align*}
for every $i\in[n]$. This concludes the  proof of Proposition~\ref{cor:M1}.  \hfill $\blacksquare$

\section{Conclusion}\label{sec:conclusions}
     We have studied non-convex distributed optimization over time-varying networks with lossy information sharing. Inspired by the original averaging-based distributed optimization algorithm, we proposed and studied a two-time scale decentralized algorithm including a damping mechanism for incoming information from the neighboring agents as well as local cost functions' gradients. We presented the convergence rate for various choices of the diminishing step-size parameters. By optimizing the achieved rate over all feasible choices for parameters, the algorithm obtains a convergence rate of $\BO({T}^{-1/3 + \epsilon})$ for non-convex distributed optimization problems over time-varying networks, for any $\epsilon>0$. Our theoretical results are corroborated by numerical simulations.
     \bibliographystyle{ieeetr}
    \bibliography{arxiv_dimix.bib}
\appendix
\section{Proof of Auxiliary Lemmas}\label{sec:proof-aux}
In this section, we provide the proofs of auxiliary lemmas.

{\it Proof of Lemma~\ref{lm:nor_squr}:} 
For two vectors $\bfu$ and $\bfv$ (with identical dimension) and any scalar $\omega>0$, we have
   \begin{align*}
       \|\bfu+\bfv\|^2 & = \|\bfu\|^2 + \|\bfv\|^2 + 2\la\bfu,\bfv \ra \cr
       &\stackrel{\rm{(a)}}{\leq}  \|\bfu\|^2 + \|\bfv\|^2 + 2 \|\bfu \| \|\bfv\| \\
       &= \|\bfu\|^2 + \|\bfv\|^2 + 2 \lp \sqrt{\omega}\|\bfu \| \cdot \frac{1}{\sqrt{\omega}}\|\bfv\|\rp \\
        &\stackrel{\rm{(b)}}{\leq}  \|\bfu\|^2 + \|\bfv\|^2 +  \omega \|\bfu \|^2 + \frac{1}{\omega}\|\bfv\|^2 \\
        &= (1+\omega)  \|\bfu\|^2 + \lp 1+ \frac{1}{\omega} \rp \|\bfv\|^2,
   \end{align*}
   where $\rm{(a)}$ follows from the Cauchy–Schwarz inequality and we used the geometric-arithmetic inequality in $\rm{(b)}$.	
   
   Now, recall that for a matrix  $U\in \mathbb{R}^{n\times d}$ we have $\|U\|_\bfr = \sum_{i=1}^n r_i \|U_i\|^2$. Therefore, for matrices $U$ and $V$ we have 
   \begin{align*}
       \|U+V\|_\bfr^2 &= \sum_{i=1}^n r_i \|U_i+V_i\|^2 \leq \sum_{i=1}^n r_i \lb (1+\omega) \|U_i\|^2 + \lp 1+\frac{1}{\omega}\rp \|V_i\|^2\rb \cr
       &= (1+\omega) \sum_{i=1}^n r_i \|U_i\|^2 + \lp 1+\frac{1}{\omega}\rp \sum_{i=1}^n r_i \|V_i\|^2\cr
       &= (1+\omega) \|U\|^2_\bfr + \lp 1+\frac{1}{\omega}\rp  \|V\|^2_\bfr.
   \end{align*}
   This completes the proof of the lemma. 
   \hfill $\blacksquare$ 
   
{\it Proof of Lemma~\ref{lemma:Psum}:} 
        To prove~\eqref{eq:esumprod}, let $g=\sum_{s=1}^{t-1} \beta(s)\prod_{k=s+1}^{t-1}(1-\lambda \beta(k))$. Then, we have
        \begin{align*}
            \lambda g &= \sum_{s=1}^{t-1} \lambda \beta(s)\prod_{k=s+1}^{t-1}(1-\lambda \beta(k))\cr
                      &= \sum_{s=1}^{t-1} (1-(1-\lambda \beta(s)))\prod_{k=s+1}^{t-1}(1-\lambda \beta(k))\cr 
                      &= \sum_{s=1}^{t-1}\left[\prod_{k=s+1}^{t-1}\lp 1-\lambda \beta(k)\rp -\prod_{k=s}^{t-1}\lp1-\lambda \beta(k))\rp\right].
        \end{align*}
        Note that the last summation is a telescopic sum, we have
        \begin{align*}
            \lambda g &= \prod_{k=t}^{t-1}\lp 1-\lambda \beta(k)\rp -\prod_{k=1}^{t-1}\lp1-\lambda \beta(k) \rp = 1 -\prod_{k=1}^{t-1}\lp1-\lambda \beta(k) \rp . 
        \end{align*}
        Dividing both sides by $\lambda > 0$ we arrive  at~\eqref{eq:esumprod}. To show \eqref{eq:esumprodb}, we utilize the same idea by letting ${h= \sum_{t=s+1}^{T} \beta(t)\prod_{k=s+1}^{t-1}(1-\lambda \beta(k))}$, we can write
        \begin{align*}
            \lambda h&=\sum_{t=s+1}^T\left[\prod_{k=s+1}^{t-1}\lp 1-\lambda \beta(k)\rp-\prod_{k=s+1}^{t}\lp 1-\lambda \beta(k)\rp\right]\cr 
            &= \prod_{k=s+1}^{s}\lp 1-\lambda \beta(k)\rp-\prod_{k=s+1}^{t}\lp 1-\lambda \beta(k)\rp\cr
            &=1-\prod_{k=s+1}^{T}\lp 1-\lambda \beta(k)\rp. 
        \end{align*}
        Dividing both sides by $\lambda>0$, results in~\eqref{eq:esumprodb}.	\hfill $\blacksquare$ 

{\it Proof of Lemma~\ref{lm:sum-exp}:} In order to prove the lemma, we separately analyze the cases $\dl<-1$, $\dl=-1$, $-1<\dl<0$, and $\dl\geq 0$. Note that for $\dl<0$,  the function  $h(t) := (t+\tau)^{\dl}$ is a decreasing function, and thus we have ${\sum_{t=1}^T (t+\tau)^\dl \leq \int_{0}^{T} (t+\tau)^{\dl}dt}$. In the following we upper bound the latter integral for each regime of $\dl$. When $\dl<-1$ we have
\begin{align}
\sum_{t=1}^T (t+\tau)^\dl \leq \int_{0}^{T} (t+\tau)^{\dl}dt = \frac{ \tau^{1+\dl} - (T+\tau)^{1+\dl}}{-1-\dl} \leq 
\frac{\tau^{1+\dl}}{-1-\dl},
\end{align}
which does not grow with $T$. 
For $\dl=-1$,  we get 
\begin{align}
\sum_{t=1}^T (t+\tau)^{-1} \leq \int_{0}^{T} (t+\tau)^{-1}dt = \ln (T+\tau) - \ln(\tau) = \ln \lp \frac{T}{\tau}+1 \rp .
\end{align}
When $-1< \dl<0 $ we arrive at 
\begin{align}
\sum_{t=1}^T (t+\tau)^\dl \leq \int_{0}^{T} (t+\tau)^{\dl}dt = \frac{(T+\tau)^{1+\dl} - \tau^{1+\dl}  }{1+\dl} \leq 
\frac{(T+\tau)^{1+\dl}}{1+\dl} \leq \frac{2^{1+\dl}}{1+\dl}(T+\tau)^{1+\dl}, 
\end{align}
where the last inequality holds since $2^{1+\dl} \geq 2^0=1$. 

Finally, for  $\dl\geq 0$, the function $h(t)$ is an increasing function. 
Hence, we can write
\begin{align}
   \sum_{t=1}^T (t+\tau)^\dl \leq \int_{1}^{T+1} (t+\tau)^{\dl}dt = \frac{(T+\tau+1)^{1+\dl} - (\tau+1)^{1+\dl} }{1+\dl}
   \leq \frac{2^{1+\dl}}{1+\dl} (T+\tau)^{1+\dl}. 
\end{align}
Note that the last inequality follows from $T+\tau+1\geq 2(T+\tau)$ which holds for any $T\geq 1$ and $\tau\geq 0$. \hfill $\blacksquare$ 

{\it Proof of Lemma~\ref{lemma:matrixprop}:} 
Recall that the $(i,j)$th entry of the matrix product $AB$ is the inner product between the $i$th row of $A$ and the $j$th column of $B$. Thus, using  the Cauchy-Schwarz inequality, we have ${|[AB]_{ij}|=|\la A_i,B^j\ra| \leq \lnr A_i\rnr \lnr B^j\rnr}$. Therefore, 
	  \begin{align*}
	   \lnr[AB]_i\rnr^2
	   =\sum_{j=1}^m |[AB]_{ij}|^2 = 
	   \sum_{j=1}^m|\la A_i,B^j\ra|^2\leq \lnr A_i\rnr^2 \sum_{j=1}^m\lnr B^j\rnr^2\leq \lnr A_i\rnr^2 \lnr B\rnr_F^2.
	  \end{align*}
	  
	 Using this inequality and the definition of $\bfr$-norm, we get
	  \begin{align*}
	     \lnr AB \rnr_\bfr^2 &=  \sum_{i=1}^nr_i\lnr[AB]_i\rnr^2\leq \sum_{i=1}^nr_i\lnr A_i\rnr^2 \lnr B\rnr_F^2=\lnr A\rnr_\bfr^2 \lnr B\rnr_F^2.
 	  \end{align*}\hfill $\blacksquare$ 
	  
{\it Proof of Lemma~\ref{lemma:transition2}:} 
Due to the separable nature of $\Nr{\cdot}$, i.e.,  $\Nr{U}^2=\sum_{j=1}^d\Nr{U^{j}}$, without loss of generality, we may assume that $d=1$. Thus, let $U=\bfu\in \R^n$. Define $\V:\R^n\to \R^+$ by 
	    \begin{align}
	        \V(\bfu) := \Nr{\bfu-\ones \bfr^T\bfu}^2=\sum_{i=1}^n r_i(u_i-\bfr^T\bfu)^2.
	    \end{align}
	    For notational simplicity, let $\bfu(s)=\bfu=\begin{bmatrix} u_1 & u_2 & \ldots & u_n \end{bmatrix}$, and  $\bfu(k+1)=A(k)\bfu(k)$. Also with a slight abuse of notation, we  denote $\V(\bfu(k))$ by $\V(k)$ for $k=s,\ldots,t$. 
	    
	    Using Theorem 1 in  \cite{touri2011existence}, we have 
	    \begin{align}\label{eqn:decreaseV}
	       \V(t)=\V(s)-\sum_{k=s}^{t-1}\sum_{i<j}H_{ij}(k)(u_i(k)-u_j(k))^2,
	    \end{align}
	    where $H(k)=A^T(k)\diag(\bfr)A(k)$. Note that $A(k)$ is a non-negative matrix, and hence we have ${H(k)\geq \rmin A^T(k)A(k)}$, for $k=s,\ldots,t$. Also, note that since ${A(k)=(1-\beta(k))I+\beta(k)W(k)}$, then  Assumption~\ref{asm:W}-(b) implies that the minimum non-zero elements of $A(k)$ are bounded bellow by $\minW\beta(k)$. Therefore, since $\beta(k)$ is non-increasing, on the window $k=s,\ldots,s+B$, the minimum non-zero elements of $A(k)$  for $k$ in this window are lower bounded by $\minW\beta(s+B)$. Without loss of generality assume that the entries of $\bfu$ are sorted, i.e., $u_1\leq \ldots\leq u_n$, otherwise, we can relabel the agents (rows and columns of $A(k)$s and $\bfu$ to achieve this).  Therefore, by Lemma~8 in \cite{nedic2008distributed}, for \eqref{eqn:decreaseV}, we have 
	    \begin{align}\label{eqn:decreaseV2}
	       \V(s+B) &\leq\V(s)- \rmin\sum_{k=s}^{s+B-1}\sum_{i<j}[A^T(k)A(k)]_{ij}(u_i(k)-u_j(k))^2\nonumber\\ 
	       &\leq \V(s)-\frac{\minW\rmin}{2}\beta(s+B)\sum_{\ell=1}^{n-1}(u_{\ell+1}-u_{\ell})^2.
	    \end{align}
	    
	    We may comment here that although Lemma~8 in \cite{nedic2008distributed} is written for doubly stochastic matrices, and its statement is about the decrease of $\V(\bfx)$ for the special case of  $\bfr=\frac{1}{n}\ones$, but in fact, it is a result on bounding $\sum_{k=s}^{s+B-1}\sum_{i<j}[A^T(k)A(k)]_{ij}(u_i(k)-u_j(k))^2$ for a sequence of {$B$-connected} stochastic matrices $A(k)$ in terms of the minimum non-zero entries of stochastic matrices ${A(s),\ldots,A(s+B-1)}$. 
	    
	    Next, we will show that $\sum_{\ell=1}^{n-1}(u_{\ell+1}-u_{\ell})^2\geq n^{-2}\V(\bfu)$. This argument adapts a similar argument used in the proof of Theorem~18 in \cite{nedic2008distributed} to the general $\V(\cdot)$. 
	    
	    For a $\bfv\in \R^n$ with $\V(\bfv)>0$, define the quotient 
	    \begin{align}
	        h(\bfv)=\frac{\sum_{\ell=1}^{n-1}(v_{\ell+1}-v_{\ell})^2}{\sum_{i=1}^nr_i(v_i-\bfr^T\bfv)^2}=\frac{\sum_{\ell=1}^{n-1}(v_{\ell+1}-v_{\ell})^2}{\V(\bfv)}. 
	    \end{align}
	    Note that $h(\bfv)$ is invariant under scaling and translations by all-one vector, i.e., $h(\omega \bfv )=h(\bfv)$ for all non-zero $\omega\in \R$ and $h(\bfv+\omega \ones)=h(\bfv)$ for all $\omega\in \R$. Therefore, 
	    \begin{align}\label{eqn:hv}
	        \min_{\substack{v_1\leq v_2\leq \cdots\leq v_n\\\V(\bfv)\not=0}}h(\bfv)&=\min_{\substack{v_1\leq v_2\leq \cdots\leq v_n\\ \bfr^T\bfv=0,\V(\bfv)=1}}h(\bfv)\cr 
	        &=\min_{\substack{v_1\leq  v_2\leq \cdots\leq v_n\\\bfr^T\bfv=0,\V(\bfv)=1}}\sum_{\ell=1}^{n-1}(v_{\ell+1}-v_{\ell})^2. 
	    \end{align}
	    
	    Since $\bfr$ is a stochastic vector, then for a vector $\bfv$ with  ${v_1\leq \ldots\leq v_n}$ and $\bfr^T\bfv=0$, we would have ${v_1\leq \bfr^T\bfv=0\leq v_n}$. On the other hand, the fact that ${\V(\bfv)=\sum_{i=1}^nr_iv^2_i=1}$ would imply ${\max(|v_1|,|v_n|)\geq \frac{1}{\sqrt{n}}}$.
	    Let us consider the difference sequence  $\hat{v}_\ell=v_{\ell+1}-v_{\ell}$  for $\ell=1,\ldots, n-1$, for which we have $\sum_{i=1}^{n-1} \hat{v}_\ell = v_n - v_1 \geq v_n \geq \frac{1}{\sqrt{n}}$. Therefore, the optimization problem \eqref{eqn:hv} can be  rewritten as
        \begin{align}\label{eqn:hv2}
	        \min_{\substack{v_1\leq v_2\leq \cdots\leq v_n\\\V(\bfv)\not=0}}h(\bfv)&=\min_{\substack{v_1\leq v_2\leq \cdots\leq v_n\\\bfr^T\bfv=0,\V(\bfv)=1}}\sum_{\ell=1}^{n-1}(v_{\ell+1}-v_{\ell})^2\cr 
	        &\geq 
	        \min_{\substack{\hat{v}_1,\ldots,\hat{v}_{n-1}\geq 0\\\sum_{i=1}^{n-1}\hat{v}_i\geq \frac{1}{\sqrt{n}}}}\sum_{\ell=1}^{n-1}\hat{v}_\ell^2. 
	    \end{align}
	    Using the Cauchy-Schwarz inequality, we get 
	    ${\sum_{\ell=1}^{n-1} \hat{v}_\ell^2 \cdot \sum_{\ell=1}^{n-1} 1^2 \geq \big( \sum_{\ell=1}^{n-1} \hat{v}_\ell \big)^2 \geq \big( \frac{1}{\sqrt{n}} \big)^2 = \frac{1}{n}}$. Hence, 
	     
	    \begin{align}
	         \min_{\substack{v_1\leq v_2\leq \ldots\leq v_n\\
	         \V(\bfv)\not=0}}h(\bfv)\geq \frac{1}{n(n-1)}\geq \frac{1}{n^2}.
	    \end{align}
        This implies that for $v_1\leq v_2\leq \ldots\leq v_n$, we have $\sum_{\ell=1}^{n-1}(v_{\ell+1}-v_{\ell})^2\geq n^{-2}\V(\bfv)$ (note that this inequality also holds for $\bfv\in\R^n$ with $\V(\bfv)=0$). 
	    Using this fact in \eqref{eqn:decreaseV2} implies  
	    \begin{align}\label{eqn:Vdot}
	       \V(s+B) \leq \lp 1-\frac{\minW\rmin}{2n^2}\beta(s+B)\rp\V(s).
	    \end{align}
	    Applying~\eqref{eqn:Vdot} for $\Del:=\lfloor \frac{t-1-s}{B}\rfloor$ steps recursively, we get
	    \begin{align*}
	        \V(s+\Del B) \leq \prod_{j=1}^{\Del} \lp 1-\frac{\minW\rmin}{2n^2}\beta(s+jB)\rp\V(s)
	    \end{align*}
	  Using the fact that  $(1-x)^{1/B} \leq 1-x/B$ and since $\{\beta(k)\}$ is a non-increasing sequence, we have
\begin{align*}
	1-\frac{\minW\rmin}{2n^2}\beta(s+jB) 
	&= \prod_{\ell=1}^{B} \lp 1-\frac{\minW\rmin}{2n^2}\beta(s+jB)\rp^{1/B}\\  
	&\leq 
	\prod_{\ell=1}^{B}  \lp 1-\frac{\minW\rmin}{2Bn^2}\beta(s+jB)\rp\nonumber\\
	&\leq 
	\prod_{\ell=1}^{B}  \lp 1-\lambda \beta(s+jB+\ell)\rp.
\end{align*}
Recall from~\eqref{eqn:decreaseV} that $\V(\cdot)$ is non-increasing. Therefore, for $s+\Del B \leq t -1 < s+(\Del+1) B $ we have
\begin{align}\label{eq:B-jumps}
	  \V(t-1) &\leq \V(s+\Del B)\nonumber\\
	  &\leq \prod_{j=1}^{\Del} \lp 
	  1-\frac{\minW\rmin}{2n^2}\beta(s+jB) 
	  \rp\V(s) \nonumber\\
	  &\leq \prod_{j=1}^{\Del} 	\prod_{\ell=1}^{B}  \lp 1- \lambda \beta(s+jB+\ell)\rp \V(s) \nonumber\\
	  &= \prod_{k=s+B+1}^{s+(\Del+1)B} \lp 1- \lambda \beta(k)\rp \V(s) \nonumber\\
	  &\leq \prod_{k=s+B+1}^{t-1} \lp 1-\lambda \beta(k)\rp \V(s). 
\end{align}
	   Next, since $\{\beta(k)\}$ is a non-increasing sequence, we have $\beta(k) \leq \beta(1) = \bzr$. Thus, 
	   \begin{align}\label{eq:initial-steps}
	       \prod_{k=s+1}^{s+B} \lp 1- \lambda \beta(k)\rp &\geq \prod_{k=s+1}^{s+B} \lp 1- \lambda \bzr \rp \nonumber\\
	       &= \lp 1- \lambda \bzr \rp^B \geq 1- B\lambda \bzr. 
	   \end{align}
	   Therefore, combining~\eqref{eq:B-jumps} and~\eqref{eq:initial-steps}, we get 
	   \begin{align}\label{eq:s-to-t}
	       \V(t-1) &\leq \prod_{k=s+B+1}^{t-1} \lp 1- \lambda \beta(k)\rp  \V(s) \nonumber\\
	       &\leq 
	       \frac{\prod_{k=s+1}^{s+B} \lp 1- \lambda \beta(k)\rp }{1-B\lambda \bzr}\!\! \prod_{k=s+B+1}^{t-1} \!\!\lp 1- \lambda \beta(k)\rp  \V(s) \nonumber\\
	       &= 
	       \frac{1}{1-B\lambda \bzr} \prod_{k=s+1}^{t-1} (1-\lambda \beta(k)) \V(s). 
	   \end{align}
	   Now, we define 
	   \[
	   \Phi(t\hspace{-1pt}:\hspace{-1pt}s)=A(t\hspace{-1pt}-\hspace{-1pt}1)\cdots A(s\hspace{-1pt}+\hspace{-1pt}1)\] for ${t\geq s}$ with ${\Phi(t\hspace{-1pt}:\hspace{-1pt}t\hspace{-1pt}-\hspace{-1pt}1)\!=\!I}$. Note that   {Assumption~\ref{asm:W}-(a)} and the fact ${A(k) = (1\!-\!\beta(k))I \!+\! \beta(k) W(k)}$ imply ${\bfr^T \Phi(t:s)  =\bfr^T}$. Then, setting $\bfu(s) = \bfu = U$ and ${\bfu(t-1) = \Phi(t:s) \bfu(s) =\Phi(t:s) U}$, we can write
	   \begin{align*}
	       \lp \Phi(t:s) -\ones \bfr^T\rp U & = \Phi(t:s)U- \ones \bfr^T U \nonumber\\
	       &=\Phi(t:s)U-\ones \bfr^T \Phi(t:s) U \nonumber\\
	       &= \bfu(t-1) - \ones \bfr^T \bfu(t-1).
	   \end{align*}
	   Therefore, using~\eqref{eq:s-to-t} we have
	   \begin{align*}
	       \big\| \big( \Phi(t:s)-\ones \bfr^T\big) U \big\|_\bfr^2  & = \Nr{\bfu(t-1) - \ones \bfr^T \bfu(t-1)}^2 = \V(t-1) \cr 
	       & \leq \frac{1}{1-B\lambda \bzr} \prod_{k=s+1}^{t-1} (1-\lambda \beta(k)) \V(s)  \\
	       &= \frac{1}{1-B\lambda \bzr} \prod_{k=s+1}^{t-1} (1-\lambda \beta(k)) \Nr{\bfu - \ones \bfr^T \bfu}^2 
	       \cr
	       & \leq \frac{1}{1-B\lambda \bzr} \prod_{k=s+1}^{t-1} (1-\lambda \beta(k)) \Nr{\bfu}^2 \\
	       &=\frac{1}{1-B\lambda \bzr} \prod_{k=s+1}^{t-1} (1-\lambda \beta(k)) \Nr{U}^2, 
	   \end{align*}
	   where the second inequality follows form the the fact that $\Nr{\bfu - \ones \bfr^T \bfu}^2 + \Nr{\ones \bfr^T \bfu}^2 = \Nr{\bfu}^2$. This completes the proof of the lemma. 
	\hfill $\blacksquare$

\end{document}